\documentclass[a4paper,12pt]{amsart}

\usepackage{amstext}
\usepackage{amsmath}
\usepackage{ifthen}
\usepackage{amscd}
\usepackage{amssymb}
\usepackage[french]{babel}
\usepackage[T1]{fontenc}
\usepackage[utf8]{inputenc}
\usepackage[all]{xy}
\usepackage{graphicx}
\usepackage{enumerate}
\usepackage{xspace}
\usepackage{pdfsync}
\usepackage{epic}
\usepackage{dsfont}
\newenvironment{demo}[1][]{\ifthenelse{\equal{#1}{}}{\noindent\textbf{Démonstration :}\xspace}{\noindent\textbf{Démonstration #1:}\xspace}}{$\square$\newline}
\newtheoremstyle{break}% name
  {}%      Space above, empty = `usual value'
  {}%      Space below
  {\itshape}% Body font
  {}%         Indent amount (empty = no indent, \parindent = para indent)
  {\bfseries}% Thm head font
  {.}%        Punctuation after thm head
  {\newline}% Space after thm head: \newline = linebreak
  {}%         Thm head spec
  
\newtheoremstyle{rq}
  {}
  {}
  {\slshape}
  {}
  {\bfseries}
  {.}
  {3pt}
  {}

\newtheoremstyle{exemple}
  {}
  {}
  {\upshape}
  {}
  {\bfseries}
  {.}
  {3pt}
  {}

\newtheoremstyle{fact}
  {}
  {}
  {\slshape}
  {}
  {\bfseries}
  {.}
  {2pt}
  {}

\theoremstyle{fact}
\newtheorem*{fact*}{Fact}
\theoremstyle{break}
\newtheorem{thmi}{Théorème}

\newtheorem{thm}{Théorème}[section]

\newtheorem{cor}[thm]{Corollaire}
\newtheorem{lem}[thm]{Lemme}
\newtheorem{prop}[thm]{Proposition}
\newtheorem{defi}[thm]{Définition}
\theoremstyle{rq}
\newtheorem{rem}[thm]{Remarque}

\theoremstyle{exemple}
\newtheorem{ex}[thm]{Exemple}

\newcommand{\VV}{\mathbb V}

\newcommand{\CC}{\mathbb C}
\newcommand{\QQ}{\mathbb Q}

\newcommand{\PP}{\mathbb P}

\newcommand{\gf}[1]{\pi_1(#1)}

\newcommand{\To}{\longrightarrow}

\newcommand{\set}[1]{\left\{#1\right\}}

\renewcommand{\ker}{\operatorname{Ker}}
\newcommand{\gl}{\operatorname{GL}}
\newcommand{\alb}{\operatorname{Alb}}

\title[Représentations linéaires des groupes kählériens]{Représentations linéaires des groupes kählériens : Factorisations et conjecture de Shafarevich linéaire}
\date{24 Janvier 2014}
\author{Fréderic Campana, Beno\^it Claudon, Philippe Eyssidieux}
\address{Fréderic \textsc{Campana}, Beno\^it \textsc{Claudon}, Université de Lorraine, Institut \'Elie Cartan Nancy, UMR 7502, B.P. 70239, 54506 Vand\oe uvre-lès-Nancy Cedex, France}
\address{Philippe \textsc{Eyssidieux}, Institut Fourier, Université Grenoble 1, 38402 Saint-Martin d'Hères Cedex, France}
\email{Frederic.Campana@univ-lorraine.fr}
\email{Benoit.Claudon@univ-lorraine.fr}
\email{Philippe.Eyssidieux@ujf-grenoble.fr}
\begin{document}

\begin{abstract}
Nous étendons aux variétés kählériennes compactes quelques 
résultats classiques sur les représentations 
linéaires des groupes fondamentaux des variétés projectives lisses. 
Notre approche, basée sur une interversion de fibrations \`a fibres tores vs
 variétés de type général, fournit une alternative \`a celle de \cite{Z96}.
Enfin nous étendons au cas kählérien les résultats généraux de convexité holomorphe pour les revêtements associés connus dans le cas projectif. 

\end{abstract}

\maketitle

\tableofcontents

\section{Introduction}

Les travaux fondamentaux de Corlette et Simpson sur les repr\'esentations lin\'eaires complexes des groupes fondamentaux
des variétés  k\"ahlériennes compactes \cite{S92,S93,S94a,S94b} ainsi que ceux de Gromov et Schoen 
sur leurs représentations à valeurs dans des corps locaux \cite{GS92} ont permis d'améliorer notre compréhension
des groupes fondamentaux des variétés algébriques lisses \cite{Z96,Z99,JZ2,LasRam} -voire aussi le survey 
\cite{ABCKT}- et des rev\^etements
associés \cite{KR,E04}. Le développement de ces idées a permis récemment d'établir 
que le rev\^etement universel d'une variété 
projective lisse complexe de groupe fondamental linéaire 
est holomorphiquement convexe \cite{EKPR} -voir aussi \cite{E10} pour un survey  r\'ecent. 

La première motivation de ce travail \'etait de généraliser au cas des variétés kählériennes compactes les résultats
de \cite{E04,EKPR}. Ceux ci sont limités au cas projectif car plusieurs résultats
sur la théorie des représentations du groupe fondamental qui sont utilisés de façon 
essentielle dans \cite{E04}
ne sont pas disponibles dans la littérature pour le cas k\"ahlérien général. 
Ces résultats sont l'ubiquité des Variations de Structure de Hodge \cite{S92}
et la théorie des ensembles
constructibles absolus de classes de conjugaison de représentations linéaires complexes \cite{S93}. 
Le présent article remédie à cette déficience de la littérature en établissant une version kählérienne
de ces outils.

Notre approche repose sur un travail important de Zuo \cite{Z96} dont
le résultat principal doit être vu comme une formulation du principe
 que la théorie des représentations linéaires complexes des groupes k\"ahlériens se 
ram\`ene au cas projectif. Nous affinons le principe de Zuo, 
en l'unifiant avec l'existence du morphisme de Shafarevich (voir la définition \ref{défi morphisme Shafarevich}):

\begin{thmi}\label{t1}
 Soit $X$ une variété k\"ahl\'erienne compacte et $\rho:\pi_1(X)\to GL_N(\mathbb{C})$ 
une représentation linéaire d'image Zariski dense dans un groupe semi-simple. Alors le morphisme de Shafarevich $sh_{\rho}: X\to Sh_{\rho}(X)$
existe et $Sh_{\rho}(X)$ est une variété projective algébrique normale de type g\'en\'eral si $\rho(\gf{X})$ est sans torsion. De plus, 
si $e:X'\to X$ est un rev\^etement étale fini tel que $e^*\rho (\gf{X'})$ est sans torsion,
$e^*\rho$ factorise par un mod\`ele lisse de $Sh_{e^*\rho}(X')$.
\end{thmi}

On notera que l'existence d'un tel $e$ est impliquée par le lemme de Selberg. 
En raison d'un point délicat dans la preuve de \cite{Z96}, nous avons   fourni
 une démonstration nouvelle de ses résultats,   basée sur \cite{E04}
 et un énoncé qui nous semble d'un intérêt indépendant:

\begin{thmi}\label{t2}
Soit $X$ une variété kählérienne 
compacte et $f:X\To K$ une application holomorphe dans une variété de Kummer
 dont les composantes connexes des fibres lisses sont  de type général. 

Il existe alors une application holomorphe génériquement finie $r:X'\To X$, une application biméromorphe
 $X'\dashrightarrow X''$ et $g:X''\To Y$ un 
fibré principal holomorphe de groupe structurel un tore complexe,  $Y$ étant une variété algébrique  de type général.

\end{thmi}

L'ingrédient principal de la preuve est l'additivité des dimensions de Kodaira 
quand la fibre est de type général dans le cas k\"ahlérien \cite{Nak}.

Ainsi formulé, le principe de Zuo permet de donner une définition d'ensemble constructible absolu
de la variété des caractères du groupe fondamental d'une variété k\"ahlérienne compacte
généralisant \cite{S93} puis de généraliser au cas kählérien les résultats 
de \cite{E04,EKPR}: 

\begin{thmi}\label{t3}
Soit $X$ une variété kählérienne compacte, $G$ un groupe algébrique linéaire 
réductif défini sur $\QQ$ et $M\subset M_B(X,G)$  un ensemble constructible absolu de la variété des caractères associée.

 Considérons $\widetilde{H^0_M}$ le sous-groupe de $\gf{X}$ défini 
comme l'intersection des noyaux de \emph{toutes} les représentations
semisimples $\gf{X}\To G(\mathbb{C})$ dont la classe de conjugaison est dans $M$ et, si $M= M_B(X,G)$, $\widetilde{H^{\infty}_M}$ le sous-groupe de $\gf{X}$ défini 
comme l'intersection des noyaux de \emph{toutes} les représentations $\gf{X}\To G(A)$ 
où $A$ est une $\CC$-algèbre arbitraire. 

Le revêtement galoisien de groupe $\widetilde{H^0_M}$ (resp. $\widetilde{H^{\infty}_M}$)
$$\widetilde{X^{*}_M}:=\widetilde{X^u}/\widetilde{H^{*}_M}, \ \ *=0,\infty$$
est alors holomorphiquement convexe.
\end{thmi}

On établira enfin un théorème de structure pour la $\Gamma$-réduction
d'une représentation linéaire d'un groupe k\"ahlérien:

\begin{thmi} \label{t4}
Soit $X$ une variété kählérienne compacte,
Soit $\rho:\gf{X}\to \gl_N(\CC)$ une représentation linéaire. 
Alors 
si $e:X'\to X$ est un rev\^etement étale fini tel que $e^*\rho (\gf{X'})$ est sans torsion,
 la $\Gamma$-réduction $Sh_{e^*\rho}(X')$ est biméromorphe 
à l'espace total d'une fibration lisse $\tau:Sh_\rho(X)\to S_\rho(X)$ 
en tores complexes sur une variété algébrique $S_\rho(X)$,  de type général.
\end{thmi}

Outre les éléments cités ci-dessus, la démonstration utilise une nouvelle fois de façon cruciale \cite{Nak}
via ses résultats
 sur les fibrations en $Q$-tores.

Décrivons maintenant le contenu de l'article. La section 2 
donne les définitions et les premières propriétés de
la $\Gamma$-réduction et du morphisme de Shafarevich. La section 3
donne une démonstration du théorème \ref{t1}. La section 4 traite les résultats 
que nous venons d'énoncer dans
le cas des représentations linéaires résolubles des groupes kählériens
comme prélude à l'étude du cas général. La section 5 établit les résultats de
convexité holomorphe du théorème \ref{t3}. La section 6 \'etablit le théorème \ref{t4}.

Nous remercions S. Druel, V. Koziarz, M. P\u{a}un et C. Voisin pour d'utiles remarques
relatives à ce travail. 

{\bf{Remarque}}: Il est naturel d'appliquer ce travail \`a une question li\'ee au probl\`eme de Serre
de caract\'eriser les groupes de pr\'esentation finie apparaissant comme groupes fondamentaux d'une vari\'ete projective complexe (groupes projectifs complexes): est ce que 
tout groupe k\"ahl\'erien (groupe fondamental d'une vari\'et\'e k\"ahl\'erienne compacte) est un groupe projectif complexe? 

La premi\`ere version diffus\'ee de ce travail arXiv:1302.5016v1 contenait une preuve erron\'ee de l'\'enonc\'e suivant qui implique que tout groupe k\"ahl\'erien
lin\'eaire est commensurable \`a un groupe projectif complexe. 

\begin{thmi}\label{kodaira}
Soit $X$ une variété kählérienne compacte et $\rho:\pi_1(X)\To \Gamma$ une représentation linéaire de son groupe fondamental. Il existe alors une variété projective lisse $Y$ et $\sigma:\pi_1(Y)\To \gl_N(\CC)$ une représentation dont l'image est un sous-groupe d'indice fini de $\Gamma$.

En particulier, si $\pi_1(X)$ est un groupe k\"ahl\'erien linéaire, il existe une
variété projective lisse ayant pour groupe fondamental un sous groupe d'indice fini de $\pi_1(X)$.
\end{thmi}

Un travail en cours de r\'edaction donnera une preuve de cet \'enonc\'e
utilisant des techniques du probl\`eme de Kodaira. Ces techniques \'etant de nature
assez diff\'erente des pr\'esentes,  nous pr\'ef\'erons  renoncer \`a traiter ici cette application. 

\section{$\Gamma$-réduction et morphisme de Shafarevich}

\subsection{$\Gamma$-réduction}

Soit $X$ une variété k\" ahlérienne compacte connexe, 
$\rho:\pi_1(X)\to \Gamma$ un morphisme de groupes surjectif, 
$u_{univ}: \tilde{X}^u \to X$ le revêtement universel 
et $u_\rho:\tilde{X}_\rho:=\ker(\rho) \backslash \tilde{X}^u \to X$ le revêtement associé \`a $\rho$.

\begin{thm}[\cite{Ca94}]\label{g-red} 
Il existe une
application presque-holomorphe propre et connexe  de variétés k\"ahlériennes
 $\tilde{g}_\rho :\tilde{X}_\rho\dasharrow \tilde Y$, $\Gamma$-équivariante, 
 telle que pour tout $\tilde{x}\in \tilde{X}_\rho$ très général  
la fibre  de $\tilde{g_{\rho}}$ passant par $\tilde{x}$
 est le plus grand sous-ensemble analytique compact et connexe de $\tilde{X}_\rho$ passant par $\tilde{x}$.
De plus, $\tilde{g}_\rho$ est
unique \`a équivalence biméromorphe près. 
\end{thm}

Rappelons qu'une application méromorphe $f:X\dashrightarrow Y$ 
est dite presque holomorphe si l'image de son lieu d'indétermination n'est pas $Y$ ; 
de façon équivalente, cela signifie qu'elle induit une fibration propre entre 
des ouverts de Zariski de $X$ et $Y$.  Rappelons également 
qu'un point d'un espace complexe
irréductible
est dit très général s'il n'est pas 
 contenu dans une  réunion dénombrable 
de sous-espaces analytiques fermés stricts.

Prenant les quotients par $\Gamma$, nous obtenons :

\begin{cor}[\cite{Ca94}, \cite{K93}]\label{gamma réduction}
Il existe une application presque holomorphe  et connexe, unique \`a \'equivalence biméromorphe près,
  $g_\rho:X\dashrightarrow G_\rho(X)$
vérifiant la propriété suivante : si $f:Z\To X$ est
 une application holomorphe (avec $Z$ un espace complexe irréductible et normal) 
dont l'image passe par un point très général de $X$, $g_\rho\circ f$ est 
constante si et seulement si la composée $\gf{Z}\stackrel{f_*}{\To}\gf{X}\stackrel{\rho}{\To}\Gamma$ est d'image finie.
\end{cor}

\begin{rem}
L'application $g_\rho$ est appelée $\Gamma$-réduction dans \cite{Ca94}
et son introduction a été motivée par l'usage des séries de Poincaré et
 le théorème de l'indice $L^2$ d'Atiyah par  \cite{Gro1}. Lorsque $X$ est projective,
 elle est introduite dans \cite{K93}
sous le nom d'application de Shafarevich associée à $\rho$. 
Nous suivrons ici  la première terminologie, réservant la terminologie de morphisme 
de Shafarevich à une notion plus précise définie  plus loin. \end{rem}

\begin{defi}
Soit $f:X\dashrightarrow Y$ une application méromorphe entre espaces complexes compacts normaux. Un modèle lisse de 
$f$ est un diagramme $X\buildrel{r}\over{\longleftarrow} \hat{X} \buildrel{\hat{f}}\over{\longrightarrow}
 \hat{Y}\buildrel{s}\over{\longrightarrow} Y$  tel que 
$r$ et $s$ sont des applications holomorphes biméromorphes propres, $\hat{X}$ et $\hat{Y}$ des 
variétés lisses et $s \circ \hat{f}=f\circ r$. 
Un modèle lisse  de $f$ est net s'il existe $p: \hat{X}\to X'$ une application holomorphe biméromorphe 
avec $X'$ lisse 
telle que tout diviseur $\hat{f}$-exceptionnel est $p$-exceptionnel.
\end{defi}

Par \cite[Lemma 1.3]{Ca04}, toute fibration méromorphe a un modèle lisse et net.

\begin{defi}
 Nous appellerons {\em géométrique} tout  quotient abélien $\Gamma$ de $\gf{X}$ vérifiant la propriété suivante: il existe
 $A$ un tore complexe et une application holomorphe $a:X\to A$
tel que $\Gamma$ apparaisse   comme l'image
de  $\gf{X}$ dans $\gf{A}$  sous $a$. 
\end{defi}

\begin{ex}\label{exab} Si $\Gamma$ est abélien, on peut donner une description simple 
de $g_\rho$ \`a l'aide du morphisme d'Albanese $\alpha_X:X\to \alb(X)$ de $X$. En effet: 
$\Gamma$ est alors un quotient de $\gf{X}^{ab}=H_1(X,\Bbb Z)$, l'abélianisé de $\gf{X}$. 
Soit $K$ le noyau du quotient $H_1(X,\mathbb{Z})\to \Gamma$, et $B\subset \alb(X)$ le plus grand sous-tore complexe $T\subset \alb(X)$ tel que $\gf{T}<\gf{\alb(X)}=H_1(X,\Bbb Z)$ soit contenu dans $K$. Soit $q:\alb(X)\to A_\rho=\alb(X)/B$ le quotient. Alors $g_\rho$ est la factorisation de Stein de la composée (non surjective, en général) $q\circ \alpha_X:X\to A_\rho$.
En particulier, si $\alb(X)$ est un tore complexe simple et si $K$ est infini, alors $G_\rho(X)$ est un point.

Le groupe  $\gf{A_\rho}=\gf{\alb(X)}/K$ est un quotient abélien géométrique de $\Gamma=\Gamma^{ab}$, 
et le noterons $\Gamma^{geom}$.

\end{ex}

\subsection{Factorisation des représentations par la $\Gamma$-réduction}

\begin{lem}\label{lemme sans torsion}
Soit $\rho:\gf{X}\To\Gamma$ un quotient de $\gf{X}$ avec $\Gamma$ sans torsion. 
Si $g_\rho:X\To Y$ désigne un modèle lisse de la $\Gamma$-réduction de $X$,
la représentation $\rho$ se factorise par $g_{\rho}$,  i.e.: il existe un morphisme de groupes rendant commutatif le diagramme:
$$\xymatrix{\gf{X}\ar[rd]_{g_{\rho*}}\ar[rr]^{\rho} && \Gamma\\
&\gf{Y}\ar[ru]&}$$
\end{lem}
Dans ces conditions, on dira aussi que $\rho$ se factorise par $G_{\rho}(X)$. 

\begin{demo}
Considérons un modèle net $g_\rho:X\To Y$ de
 la $\Gamma$-réduction de $X$ et munissons $Y$ de la structure orbifolde 
(pour les multiplicités classiques) induite par $\Delta$ \cite{Ca04,Ca11j,Ca11}. L'adjonction du $\QQ$-diviseur $\Delta$ a pour effet de rendre la suite des groupes fondamentaux exacte :
$$\gf{X_y}\To \gf{X}\To\gf{Y/\Delta}\To 1.$$
D'autre part, le groupe fondamental orbifolde de $(Y/\Delta)$ se présente naturellement sous la forme :
$$1\To T\To \gf{Y/\Delta}\To\gf{Y}\To 1$$
où $T$ est un groupe engendré par des éléments de torsion. Par définition de $g_\rho$, $\rho(\gf{X_y})$ 
est fini donc trivial et $\rho$ se factorise donc par $\gf{Y/\Delta}$. A nouveau, comme $\Gamma$ est sans torsion, l'image $T$ par $\rho$ doit être triviale,
 ce qui signifie exactement que $\rho$ se factorise par $\gf{Y}$.
\end{demo}

Le lemme de Selberg stipule que  si $\Gamma< \mathrm{Gl}_N(\Bbb C)$ est de type fini,
 $\Gamma$ admet un sous-groupe d'indice fini $\Gamma'$ sans torsion.
Combiné avec le lemme \ref{lemme sans torsion}, il permet de déduire:

\begin{cor}\label{corselberg}
Soit $\rho:\gf{X}\To\Gamma<\mathrm{Gl}_N(\Bbb C)$ un quotient linéaire de $\gf{X}$.
 Il existe un revêtement étale fini $e:X'\to X$
 tel que si $\rho':\pi_1(X')\to \Gamma':=\rho(e_*\pi_1(X'))$ est la restriction de $\rho$, 
alors   $\rho'$  se factorise par un modèle lisse $g_{\rho'}:X'\to Y':=G_{\rho'}(X')$ de la 
$\Gamma'$ réduction de $X'$.
 Cette propriété subsiste pour tout revêtement étale $X''\to X$ dominant $X'$.\end{cor}

\subsection{Fonctorialité, Groupe quotient} 

Les propriétés élémentaires suivantes de fonctorialité et de représentation quotient 
se déduisent sans difficulté de la propriété de base de la $\Gamma$-réduction
et nous en ommettrons la preuve. 

Introduisons  les notations suivantes : si $f:U\to V$ est une application 
holomorphe propre (non nécessairement surjective) entre espaces analytiques complexes, avec $U$ normal et connexe, nous noterons $St(f):U\to St(U/V)$ et $st(f):St(U/V)\to V$ sa factorisation de Stein où $St(U/V)$ est normal, $f=St(f)\circ st(f)$, $St(f)$ est à fibres connexes et $st(f)$ est finie, d'image $f(X)$.

\begin{lem}\label{gred-imrec} 
Soit $f:W\to X$ une application holomorphe surjective entre variétés compactes K\" ahlériennes connexes, 
et $\rho: \gf{X}\to \Gamma$ un quotient du groupe fondamental de $X$.
 Soit $\rho_f:=f^*(\rho)=\rho\circ f_*:\gf{W}\dashrightarrow \Gamma_f:=Im(\rho_f)$. 
 Alors $g_{\rho_f}:W\dashrightarrow G_{\rho_f}(W)$ est la factorisation de Stein de 
la composée $g_\rho\circ f:W\dashrightarrow G_\rho(X)$ i.e.: $g_{\rho_f}=St(g_\rho\circ f)$.
\end{lem}

\begin{lem}\label{gred-comp} 
Soit $\rho:\pi_1(X)\to \Gamma<\mathrm{Gl}_N(\Bbb C)$ un quotient linéaire de $\gf{X}$. 
Soit $g_\rho:X\dashrightarrow Y:=G_\rho(X)$ sa $\Gamma$-réduction. 
Soit $\Delta<\Gamma$ un sous-groupe normal 
et $\bar \Delta$ son adhérence de Zariski dans $\mathrm{Gl}_N(\Bbb C)$. 
Notons $H$ le sous-groupe normal  de $\Gamma$ défini par $H:=\Gamma\cap \bar \Delta$.
La représentation  composée  $\sigma: \pi_1(X)\to  \Gamma/H<\bar \Gamma/\bar \Delta$ 
est alors un quotient linéaire de $\gf{X}$. Soit
 $g_\sigma:X \dashrightarrow  Z:=G_\sigma(X)$ la $\Gamma/H$-réduction de $X$. 
\begin{enumerate}[(1)]
\item Il existe alors une application méromorphe dominante $g_{\rho/\sigma}: Y:=G_\rho(X)\dashrightarrow Z:=G_\sigma(X)$ telle que $g_{\sigma}= g_{\rho/\sigma}\circ g_{\rho}$ et la restriction $g_{\rho/\sigma,z}:X_z\to Y_z$ à une fibre très générale $X_z$ de $g_\sigma$ n'est autre que la $\Gamma_z$-réduction de $X_z$ où $\sigma_z$ désigne la restriction  de $\sigma$ à $\pi_1(X_z)$ et $\Gamma_z=\sigma_z(\pi_1(X_z))$.
\item Si les images de $\rho$ et $\sigma$ sont sans torsion, ce qui est le cas quitte à remplacer $X$ par un revêtement étale fini,
 $\rho$, $\sigma$ et $\sigma_z$ 
se factorisent respectivement par $Y$, $Z$, et $Y_z$.
\item Si $\rho$ (et donc $\sigma$) se factorise par $g_\rho$ 
et si $\rho^*:\pi_1(G_\rho(X))\to \Gamma$ (resp. $\sigma^*:\pi_1(G_\rho(X))\to \Gamma/H$) factorisent $\rho$ (resp. $\sigma$), alors: $G_{\sigma^*}(G_\rho(X))=G_\sigma(X)$ et $g_{\rho/\sigma}=g_{\sigma^*}$.
\end{enumerate}
\end{lem}

Le diagramme commutatif correspondant est:

$$\xymatrix{X\ar[rd]_{g_{\sigma}}\ar[rr]^{g_\rho} && G_\rho(X):=Y\ar[ld]^{g_{\rho/\sigma}}\\
& G_\sigma(X):=Z&}$$

Nous utiliserons ce lemme 
 de différentes manières: en prenant pour $\bar \Delta$ le radical
 résoluble de $\bar \Gamma$,  $\bar \Gamma/\bar \Delta$ étant  alors semi-simple,
mais aussi en considérant 
pour $\Delta$ l'image par $\rho$ de $\pi_1(X_w)$, 
si $X_w$ est la fibre générale d'une fibration presque-holomorphe $f:X\to W$ (sans lien avec $\rho$, a priori).

\subsection{Quotient par un groupe abélien}

Nous allons donner une description de la situation du lemme \ref{gred-comp} lorsque $H$ est abélien. C'est aussi une version relative de l'exemple \ref{exab} ci-dessus.

Rappelons (\emph{cf.} \cite{Ca85}) que si $f:X\to S$ est une application holomorphe surjective à fibres connexes avec $X$ k\"ahlérienne compacte, il existe une application d'Albanese relative $\alpha_{X/S}:X\dasharrow \mathrm{Alb}(X/S)$,
 au-dessus de $S$, presque holomorphe au-dessus du lieu de lissité de $f$, dans laquelle : $\mathrm{Alb}(X/S)$ est kählérienne compacte, $\alpha(f):\mathrm{Alb}(X/S)\to S$ est holomorphe à fibres connexes avec $f=\alpha(f)\circ \alpha_{X/S}$, telle que, si $X_s:=f^{-1}(s), s\in S$ est lisse, alors $\mathrm{Alb}(X/S)_s=\alpha(f)^{-1}(s)$ est isomorphe à $\mathrm{Alb}(X_s)$ (variété d'Albanese de $X_s$) et telle que la restriction de $\alpha_{X/S}$ à $X_s$ est un morphisme d'Albanese pour $X_s$. 
En général, $\alpha_{X/S}$ n'est ni surjective ni connexe. 

\begin{lem}\label{lem ss grpe abelien} Soit $f:X\to Z$ une application holomorphe surjective \`a fibres connexes,
$z\in Z$ un point très général 
et $H$ un quotient géométrique  de $\pi_1(X_z)^{ab}$.
Il existe alors une application presque-holomorphe surjective
$$q_H:\mathrm{Alb}(X/Z)\to A_H(X/Z)$$
au-dessus de $Z$, dans laquelle $\tau_H:A_H(X/Z)\dashrightarrow Z$ est une fibration dont les fibres lisses sont des tores, 
et telles que la restriction 
$$q_{H,z}:\mathrm{Alb}(X/Z)_z \dashrightarrow\ A_H(X/Z)_z$$
de $q_H$ au-dessus de $z$ est le quotient du tore $\mathrm{Alb}(X/Z)_z$ sur le tore $A_H(X/Z)_z$ dont le groupe fondamental est $H$, quotient de $\gf{\mathrm{Alb}(X/Z)_z}$. 
\end{lem}

\begin{demo} La donnée de $H$ détermine \`a translation près un unique sous-tore $B_z$ de $\mathrm{Alb}(X/Z)_z$, 
variant holomorphiquement avec $z\in Z$ général. La projection sur son espace de paramètres du graphe 
de la famille universelle (dans l'espace des cycles de Chow-Barlet de $\mathrm{Alb}(X)$) de sous-tores relatifs 
de $\mathrm{Alb}(X/Z)$ sur $Z$ est l'application $q_H$, l'espace de paramètres étant $A_H(X/Z)$.
\end{demo}

Comme conséquence immédiate du lemme \ref{gred-comp}, nous obtenons l'énoncé suivant.

\begin{lem}\label{gred-ab} Soit $\rho:\gf{X}\to \Gamma$ un quotient (avec $X$ kählérienne compacte et connexe). Soit $H\lhd\Gamma$ un sous-groupe abélien normal et $\sigma:=s\circ \rho:\gf{X}\to \Gamma/H$ le quotient correspondant. Nous avons  un diagramme commutatif:
$$\xymatrix{X\ar[rd]_{g_{\sigma}}\ar[rr]^{g_\rho} && G_\rho(X):=Y\ar[ld]^{g_{\rho/\sigma}}\\
& G_\sigma(X):=Z&}$$

où $g_{\rho}=St(q_{H^{geom}}\circ \alpha^*_{X/Z})$ et où $q_{H^{geom}}:\mathrm{Alb}(X/Z)\to A_{H^{geom}}(X/Z)$ étant le quotient naturel associé au groupe $H^{geom}$ par le lemme \ref{lem ss grpe abelien}.
\end{lem}

\subsection{Morphisme de Shafarevich}

Si $\widetilde{X_\rho}$ est holomorphiquement convexe, il admet une réduction de Remmert :
$r:\widetilde{X_\rho}\To R\left(\widetilde{X_\rho}\right),$
c'est à dire une application $\Gamma$-équivariante propre sur un espace (normal) de Stein. Le morphisme quotient $g_\rho:X\To (R\left(\widetilde{X_\rho}\right)/\Gamma)$
est alors une $\Gamma$-réduction possédant la propriété suivante:

\begin{defi}\label{défi morphisme Shafarevich}
 Un modèle $g_\rho:X\To Y$, $Y$ pouvant être un espace complexe normal,   de la $\Gamma$-réduction sera appelé \emph{morphisme de Shafarevich} associé à $\rho$ si:
\begin{enumerate}[(i)]
\item $g_\rho$ est holomorphe,
\item pour toute application holomorphe $f:Z\To X$ (avec $Z$ espace complexe compact connexe), $g_{\rho}\circ f$ est un point si et seulement si l'image de $\gf{Z}$ par $\rho\circ f_*$ est finie.
\end{enumerate}

On notera cet unique modèle: $sh_\rho:X\to Sh_\rho(X)$ lorsqu'il existe.
\end{defi}

\begin{rem} L'existence d'un morphisme de Shafarevich est un propriété strictement plus faible 
que la convexité holomorphe de $\widetilde{X_\rho}$: il suffit qu'existe une fibration holomorphe
propre $\widetilde{X_\rho}\to R(\widetilde{X_\rho})$ telle que 
$R(\widetilde{X_\rho})$ ne contienne pas de sous-espace analytique compact
 de dimension strictement positive (sans être pour autant Stein). 
\end{rem}

\begin{ex}\label{exab'} Lorsque $\Gamma$ est abélien, il existe toujours un morphisme de Shafarevich: celui donné dans l'exemple \ref{exab}. Par contre, $\widetilde{X_\rho}$ n'est pas toujours holomorphiquement convexe dans ce cas (il existe des quotients sans fonction holomorphe non constante de $\Bbb C^n, n\geq 2$ par des réseaux partiels, dits \og groupes de Cousin\fg). 
\end{ex}

\begin{ex}\label{exab"} \label{convhol} Lorsque $\Gamma$ est abélien géométrique, 
alors $\widetilde{X_\rho}$ est holomorphiquement convexe.
\end{ex}

\section{Réduction  au cas projectif de l'étude des représentations linéaires semi-simples des groupes kählériens}

Cette section  donne une preuve alternative 
d'un  théorème fondamental de factorisation dû à Zuo \cite{Z96} (théorème \ref{factorisation presque simple}) et 
en dérive quelques conséquences connues mais non documentées dans la littérature 
comme le théorème d'ubiquité de Simpson
dans le cas k\"ahlérien.
Notre démonstration du théorème 
 necessite plusieurs étapes.
\begin{enumerate}[1.]
\item Cas d'une représentation $\rho$ associée à une variation de structures de Hodge complexes, lorsque la monodromie est discrète. Alors $Sh_\rho(X)$ est de type général.
\item Cas d'une représentation semi-simple. Nous utiliserons ici l'énoncé d'interversion des fibrations du théorème \ref{échange fibration}.
\item Cas d'une représentation réductive.
\end{enumerate}

\subsection{Interversion des fibrations}
Nous montrons ici un résultat qui peut être vu comme un phénomène d'échange des 
rôles fibre/base dans une fibration et susceptible d'applications autres que celle donnée dans le présent article.

\begin{thm}\label{échange fibration}
Soit $X$ une variété kählérienne 
compacte et $f:X\To K$ une application holomorphe dans une variété de Kummer
 dont les composantes connexes des fibres lisses sont  de type général. 

Il existe alors une application holomorphe génériquement finie $r:X'\To X$, une application biméromorphe
 $X'\dashrightarrow X''$ et $g:X''\To Y$ un 
fibré principal holomorphe de groupe structurel un tore complexe,  $Y$ étant une variété algébrique  de type général.
\end{thm}
Par variété de Kummer, on entend un quotient d'un tore complexe compact 
par un groupe fini d'automorphismes.
  
\begin{rem}\label{rq échange}
Comme on pourra l'observer dans la démonstration ci-dessous, 
la fibration $X'\To Y$ n'est autre que la fibration d'Iitaka-Moishezon de $X'$. 
Il est bien connu que la base de la fibration correspondante pour $X$ n'est pas toujours de type général. 
En effet, il suffit de considérer $X=C\times E/\langle \sigma,\tau\rangle$ où $C$ 
est une courbe hyperelliptique (de genre $g\ge2$), $\sigma$ l'involution correspondante et $E$ une courbe elliptique munie d'un point $\tau$ d'ordre 2. La deuxième projection donne (en passant au quotient) une fibration sur une courbe elliptique (en courbe hyperbolique) :
$$X\To E/\langle\tau\rangle$$
alors que la première projection fournit la fibration d'Iitaka-Moishezon :
$$X\To C/\langle\sigma\rangle\simeq\PP^1.$$

Le cas o\`u $X$ est projective est sans intérêt et résulte de l'existence d'un revêtement ramifié de $X$ qui est
 de type général. Il est standard  \cite{U75}  que, dans le théorème \ref{rq échange},
 $X$ et $X'$ ont la même
dimension algébrique $a(X)$ et, quitte \`a faire un revêtement fini supplémentaire,  on peut encore
assurer que $\dim(Y)=a(X)$ et que $g$ est la réduction algébrique de $X'$.

On aimerait, sous des hypoth\`eses suppl\'ementaires sur $f$, pouvoir conclure que $r$ peut \^etre pris \'etale. 
Nous n'avons pas trouv\'e 
de formulation raisonnable et n'obtiendrons cette pr\'ecision que dans les cas tr\`es particuliers o\`u
ce th\'eor\`eme sera appliqu\'e.
\end{rem}

\begin{demo}[du théorème \ref{échange fibration}]
Comme $K$ est une variété de Kummer, nous pouvons d'ores et 
déjà opérer un changement de base (génériquement fini) 
pour nous ramener à une application $f:X\To T$ où $T$ est 
un tore et où les composantes connexes des fibres générales sont de type général.
L'image de $X$ dans $T$ est une sous-variété irréductible $Z:=f(X)$ 
de dimension de Kodaira $\kappa(Z)\ge 0$. L'additivité des dimensions de Kodaira 
pour les fibres de type général (\cite{Ko} dans le cas projectif, 
 \cite[Th. 5.7]{Nak} dans le cas K\"ahlérien)
 montre alors que :
$$\kappa(X)\ge \kappa(X_t)+\kappa(Z)=\dim(X_t)+\kappa(Z)>0$$
où $X_t$ désigne la fibre générale de $f$ (supposée de dimension strictement positive). 
Cela signifie que la fibration d'Iitaka-Moishezon de $X$ est non constante ; 
notons la $J:X\To Y$. Si $X_y$ désigne la fibre  de $J$ en $y\in Y$ général, elle vérifie $\kappa(X_y)=0$ 
et son image dans $T$ est  un translaté d'un sous-tore $A_y$ de $T$. 
En effet, $X_y$ étant spéciale \cite{Ca11j}, son image dans $T$ l'est également et le théorème d'Uneo \cite[Th. 10.9]{U75}
 montre que les seules sous-variétés spéciales d'un tore sont les translatés de sous-tores. 
Par rigidit\'e des sous-tores, le sous-tore $A_y$ est indépendant de $y$ et on note désormais $A=A_y$.

Il est clair que $Z$ est invariant par l'action de $A$ par translation sur $T$.
Considérons alors $p_A:T\To B:=T/A$ le quotient de 
$T$ par $A$ et $W\subset B$ l'image de $X$ par $p_A\circ f$. 
Les fibres de $J$ étant envoyées sur des points par $p_A\circ f$, 
il existe une application $Y\To W$ rendant le diagramme suivant
$$\xymatrix{X\ar[d]_J\ar[r]^{f} & Z\ar[d]^{p_A}\\
Y\ar[r] & W
}$$
commutatif et nous pouvons donc examiner l'application $\mu:X\To Y\times_W Z$ 
(qui est surjective par définition de $A$). Nous allons montrer que $\mu$ est génériquement finie.
 Cela revient à montrer que $f$ est génériquement finie en restriction à $f_y:X_y\To A_y:=f(X_y)$ ; 
supposons que cela ne soit pas le cas et considérons $\tilde{f_y}:=St(f_y):X_y\To \tilde{A_y}:=St(X_y/A_y)$ la factorisation de Stein de $f_y$. Comme $y$ est général, les fibres (lisses) de $\tilde{f_y}$ sont des sous-variétés générales des fibres (générales) de $f$ : elles sont donc elles aussi de type général. Nous pouvons à nouveau appliquer l'additivité :
$$0=\kappa(X_y)\ge \kappa(\tilde{F_y})+\kappa(\tilde{A_y})=\dim(\tilde{F_y})$$
(où $\tilde{F_y}$ désigne la fibre générale de $\tilde{f_y}$).
 Nous obtenons la contradiction souhait\'ee si $f_y$ n'est pas génériquement finie.
 En particulier, $X_y$ est biméromorphe à $\tilde{A_y}$ qui est un revêtement
 étale fini de $A_y$ par
 \cite[Theorem 22, p269]{K81} (voir aussi \cite[Prop. 5.3]{Ca04}).

Le sous-groupe $\gf{\tilde{A_y}}$ est donc d'indice fini dans $\gf{A}$, lui même contenu dans $\gf{T}$ avec $\gf{B}$ pour quotient. Il existe alors un sous-groupe d'indice fini $\Gamma'\le\gf{T}$ dont l'intersection avec $\gf{A}$ coïncide avec $\gf{\tilde{A_y}}$ et vérifiant $\Gamma'/\gf{A}=\gf{B}$. Considérons le revêtement $\pi:T'\To T$ correspondant au sous-groupe $\Gamma'$ et $p:X'\To X$ le revêtement étale obtenu par changement de base. La factorisation de Stein de la composée $X'\stackrel{p}{\To} X\stackrel{J}{\To}Y$ permet de compléter le diagramme :
$$\xymatrix{X'\ar[d]_{J'}\ar[r]^{f'} & Z'\ar[d]^{\pi\circ p_A}\\
Y'\ar[r] & W
}$$
Dans cette nouvelle configuration, l'application $\mu':X'\To Y'\times_{W} Z'$ est maintenant biméromorphe et la fibration $J':X'\To Y'$ est ainsi obtenue comme image réciproque de la fibration en tore $Z'\To W$ (de fibre $\tilde{A}$) par le changement de base $Y'\To W$. En particulier, le fibré canonique relatif $K_{X'/Y'}$ est trivial et cela implique
$$\dim(Y')=\dim(Y)=\kappa(X)=\kappa(X')=\kappa\left(X',(J')^*K_{Y'}\right)=\kappa(Y'),$$
c'est à dire que $Y'$ est bien de type général.
\end{demo}

\subsection{Variétés des caractères et correspondance de Simpson}\label{section Simpson}

Rappelons quelques faits de base sur l'espace de modules des représentations linéaires d'un groupe de type fini.
 Si $\Gamma$ désigne un groupe de type fini et $G$ un groupe algébrique linéaire défine sur $\bar{\QQ}$, 
on notera $R_B(\Gamma,G)$ le schéma affine représentant le foncteur $S\mapsto \mathrm{Hom}(\Gamma,G(S))$. 
Le sch\'ema des caractères $M_B(\Gamma,G)$ est défini comme le quotient \textsc{git} $R_B(\Gamma,G)//G$ 
($G$ agissant par conjugaison) ; il s'agit donc d'un schéma affine dont les points sur un
 corps algébriquement clos de caractéristique nulle
 $\bar{k}$ sont représentés par les classes de conjugaison des représentations réductives (voir, par exemple, \cite{LM}). 
Tout point de $M_B(\Gamma,G)(\bar{k})$ sera systématiquement représenté par une telle classe de conjugaison.

Si $\Gamma=\gf{X,x}$, nous utiliserons la notation transparente $M_B(X,G):=M_B(\Gamma,G)$ (le changement de point base s'effectuant \emph{via} un automorphisme intérieur, il ne joue aucun rôle dans la définition de $M_B$).

\begin{rem}\label{remarque variété caractère}
Soit $G$ et $G'$ des groupes réductifs, $X$ et $Y$ des variétés kählériennes compactes. 
Toute application méromorphe $f:X\dashrightarrow Y$ induit un morphisme entre les schémas de caractères $f^*:M_B(Y,G)\To M_B(X,G)$. De même, tout morphisme $i:G\To G'$ induit un morphisme $i_*:M_B(X,G)\To M_B(X,G')$.
\end{rem}

Pour terminer, nous revenons sur la correspondance de Simpson
 et précisons les résultats disponibles dans la catégorie kählérienne. 
Cette correspondance établit une équivalence de catégories entre les 
représentations réductives du groupe fondamental et les fibrés de Higgs polystables
 à classes de Chern nulles (voir \cite{S92,S94b} pour les notions utilisées).

\begin{thm}\label{correspondance Simpson}
Soit $(X,\omega)$ une variété kählérienne compacte et $(E,\theta)$ un fibré de Higgs polystable vérifiant
$$\int_X c_1(E)\wedge\omega^{n-1}=\int_X c_2(E)\wedge\omega^{n-2}=0.$$
Le fibré $E$ provient alors d'une représentation réductive du groupe fondamental de $X$.
\end{thm}

Un des grands succès de la théorie de C. Simpson réside en la construction pour $X$ projective  d'un espace de modules 
$M_{Dol}(X,G)$ de classes d'isomorphisme de $G$-fibrés de Higgs polystables 
(avec $c_1=c_2=0$) \cite{S94a,S94b}. Or, cet espace porte une action naturelle\footnote{Ce phénomène a été découvert par \cite{Hit} en dimension 1.} 
de $\CC^*$ ; si $[(E,\theta)]\in M_{Dol}(X,G)$ et $t\in\CC^*$, on pose :
$$t\cdot [(E,\theta)]:=[(E,t\theta)].$$
De plus, la correspondance de Simpson s'incarne en un homéomorphisme des espaces topologiques sous-jacents :
$$M_{Dol}(X,G)(\CC)\stackrel{\sim}{\To}M_B(X,G)(\mathbb{C}).$$

Si $X$ est seulement supposée kählérienne, l'espace $M_{Dol}(X,G)$ n'a pas d'existence \emph{a priori} mais l'action de $\CC^*$ persiste
\footnote{Prendre garde au fait que cette action n'est pas algébrique en général.} sur $M_B(X,G)$ : si $[\rho]\in M_B(X,G)$ et $(E,\theta)$ est le fibré provenant de $\rho$, nous noterons $[\rho_t]$ la (classe d'isomorphisme de la) représentation associée à $(E,t\theta)$. Nous appellerons $([\rho_t])_{t\in\CC^*}$ la famille des déformations de Simpson de $[\rho]$.

Si $X$ est projective, l'action de $\CC^*$ sur la variété quasi-projective $M_{Dol}(X,G)$ est algébrique et la limite
$$\lim_{t\to 0} [(E,t\theta)]$$
existe dans $M_{Dol}(X,G)$. La classe d'isomorphisme correspondante est alors un point fixe de 
l'action de $\CC^*$ : il s'agit de la classe de conjugaison d'une représentation
 sous-jacente à une variation de structures de Hodge polarisables 
(nous utiliserons l'acronyme $\CC$-\textsc{vsh} dans la suite 
en dépit de son inélégance) \cite{S92}. 
Ce phénomène est connu sous le nom d'ubiquité des variations de structures de Hodge : 
toute représentation du groupe fondamental de $X$ peut être déformée en une  $\CC$-\textsc{vsh}.
En effet toute repr\'esentation lin\'eaire sur $\mathbb{C}$
se d\'eforme \`a sa semisimplifi\'ee  qui est réductive.

\subsection{Factorisation : cas d'une $\CC$-\textsc{vsh} discrète}
Nous abordons maintenant le problème de factorisation des représentations linéaires des groupes kählériens. Comme annoncé ci-dessus, la première étape consiste à examiner le cas particulier d'une $\CC$-\textsc{vsh} discrète. Soit donc $(\VV,\mathcal{F}^\bullet, S)$ une $\CC$-\textsc{vsh} polarisée sur une variété kählérienne compacte $X$ ; nous noterons $\rho:\gf{X,x}\To U(\VV_x,S_x)$ 
la représentation de monodromie associée au choix d'un point base $x\in X$.

\begin{prop}\label{cas discret}
Dans la situation décrite ci-dessus, supposons de plus que
 la représentation de monodromie est \emph{discrète}. 
Le morphisme de Shafarevich $X\To Sh_\rho(X)$ associé \`a $\rho$ existe  et $Sh_\rho(X)$ est de type général à revêtement étale fini près. Plus précisément, si $e:\bar{X}\To X$ est un revêtement étale fini tel que $\rho(\gf{\bar{X}})$ soit sans torsion, alors $Sh_{e^*(\rho)}(\bar{X})$ est de type général.
\end{prop}

\begin{demo}
Remarquons tout d'abord que nous pouvons supposer que l'image de $\gf{X}$ par $\rho$ est sans torsion. Considérons alors l'application de période de $(\VV,\mathcal{F}^\bullet, S)$. Il s'agit d'une application holomorphe $\rho$-équivariante :
$$\tilde{X}\To \mathcal{D}=U(\VV_x,S_x)/V$$
avec $V=\prod_k U(\VV_x^{k,-k})$. L'hypothèse sur le caractère discret de $\rho$ montre que cette application descend en un morphisme $X\To \mathcal{D}/\rho(\gf{X})$ dont la factorisation de Stein n'est autre que le morphisme de Shafarevich relativement à $\rho$ :
$$sh_\rho:X\To Sh_\rho(X),$$
$Sh_\rho(X)$ étant alors une variété projective normale (consulter \cite[p. 524-524]{E04}) que nous désignerons par $Y$ pour ne pas alourdir les notations.

La $\CC$-\textsc{vsh} $(\VV,\mathcal{F}^\bullet, S)$ descend\footnote{par $\CC$-\textsc{vsh} sur une variété normale, nous entendons un système local polarisé muni d'une filtration holomorphe vérifiant la décomposition de Hodge et la transversalité de Griffiths sur le lieu lisse.} sur $Y$ en une $\CC$-\textsc{vsh} notée $(\VV_Y,\mathcal{F}^\bullet, S)$. Si $\pi:Y^*\To Y$ est une désingularisation de $Y$, $\pi^*(\VV_Y,\mathcal{F}^\bullet, S)$ est une $\CC$-\textsc{vsh} dont l'application de période $p$ est génériquement immersive. Le fibré équivariant $T_\mathcal{D}$ descend également sur $Y^*$ en un fibré $E$, la différentielle de l'application de période devenant un morphisme de fibrés $p_*:T_{Y^*}\To E$. Remarquons que, quitte à remplacer $Y^*$ par un autre modèle birationnel, nous pouvons supposer que l'image de $T_{Y^*}\To E$ est contenu dans un sous-fibré $F$ de $E$ et que $p_*$ est génériquement un isomorphisme (au dessus d'un ouvert $U$).

Cependant, la métrique de Hodge induit une métrique $h$ sur le fibré $F$ et les formules de Griffiths et Schmid pour la courbure de $\mathcal{D}$ montrent que les courbures bissectionnelles holomorphes sont semi-négatives en restriction à $p_*T_{Y^*}$ \cite[Corollaire 9.2.2]{E04}. La courbure décroissant dans les sous-fibrés, nous en déduisons que $(F,h)$ est semi-négatif au sens de Griffiths ; en particulier $\mathrm{Tr}(i\Theta_h(F))\le0$. D'autre part, la courbure sectionnelle holomorphe étant négative dans les directions horizontales, pour tout $y\in U$ et $\alpha\in T_{Y^*,y}$ ($\alpha\neq0$), le vecteur $v=p_*(\alpha)$ satisfait $i\Theta_{\alpha\bar{\alpha}v\bar{v}}<0$ et donc $\mathrm{Tr}(i\Theta_h(F))_y<0$.

En particulier, le fibré en droites $\mathrm{det}(F^*)$ est \emph{big} (et \emph{nef}). Or, par construction, le fibré canonique de $Y^*$ se décompose en $K_{Y^*}=\mathrm{det}(F^*)+D$ où $D$ est un diviseur de Cartier effectif ; il s'ensuit que $Y^*$ est bien de type général.
\end{demo}

\begin{rem}
On pourra consulter également \cite{BKT} où des arguments  similaires 
sont utilisés pour la construction de différentielles symétriques holomorphes non triviales.
\end{rem}

\subsection{Cas semi-simple}

Dans le cas des représentations linéaires réductives définies sur un corps local, Katzarkov et Zuo 
ont développé de façon indépendante \cite{Kat95,Z96,Z99} 
une notion de revêtement spectral en utilisant la théorie des applications harmoniques à valeurs dans les 
immeubles de Bruhat-Tits \cite{GS92}. L'exploitation de ces idées est poussée plus avant dans \cite{E04}, voir également \cite{E10}.
Nous résumons les résultats de \cite[Prop. 3.4.15, Lem. 4.2.3]{E10} dont nous ferons usage sous la forme suivante.

\begin{lem}\label{factorisation KZ}
Soit $X$ une variété kählérienne compacte, $L$ un corps de nombres et $\wp$ un idéal premier de l'anneau des entiers de $L$. Si $\rho:\gf{X}\To \mathrm{GL}_N(L_\wp)$ désigne une représentation réductive, il existe une fibration holomorphe
$$s_\rho:X\To S_\rho(X)$$
vérifiant les propriétés suivantes :
\begin{enumerate}[(1)]
\item $S_\rho(X)$ est une espace kählérien normal (projectif si $X$ l'est),
\item si $Z\subset X$ désigne un sous-espace connexe de $X$, $s_\rho(Z)$ est un point si et seulement si $\rho(\gf{Z})$ est contenu dans un sous-groupe compact de $\mathrm{GL}_N(L_\wp)$.
\end{enumerate}
Si $T$ désigne une variété algébrique irréductible définie sur $\bar{\QQ}$ et si $r:T\dashrightarrow M_B(X,\mathrm{GL}_N)$ est une application rationnelle (définie elle aussi sur $\bar{\QQ}$), considérons $\rho_T:\gf{X}\To \mathrm{GL}_N(\bar{\QQ}(T))$ une représentation réductive définie sur le corps des fractions de $T$ dont la classe de conjugaison corresponde au point générique de l'image de $r$. Il existe alors une fibration holomorphe
$$s_T:X\To S_T(X)$$
qui satisfait de plus :
\begin{enumerate}[(i)]
\item $S_T(X)$ est une espace kählérien normal (projectif si $X$ l'est),
\item si $Z\subset X$ désigne un sous-espace connexe 
de $X$, $s_T(Z)$ est un point si et seulement si l'application $T\dashrightarrow M_B(Z,\mathrm{GL}_N)$ 
est constante.
\end{enumerate}
\end{lem}

\begin{rem}\label{rem KZ}
Il est à noter que les constructions effectuées dans \cite{E04} 
montrent que les fibrations $s_\rho$ et $s_T$ sont obtenues comme factorisation de Stein d'une application vers une variété de Kummer.
\end{rem}

Nous pouvons maintenant finir notre preuve alternative du résultat de factorisation de Zuo 
\cite{Z96}. 

\begin{thm}\label{factorisation presque simple}
Soit $X$ une variété kählérienne compacte et $\rho:\gf{X}\To S$
 une représentation Zariski dense dans un groupe semi simple $S$
 (défini sur un corps algébriquement clos de caractéristique nulle). 
 Il existe $\pi:X'\To X$ composition d'une modification propre et d'un revêtement étale fini,
 $f:X'\To Y$ une fibration sur une variété projective algébrique et $\rho_Y:\gf{Y}\To S$ une représentation (Zariski dense) telle que $\pi^*\rho=f^*\rho_Y$.
\end{thm}
Ceci implique que la base de la $\Gamma$-réduction de $X$ est  algébrique. \cite{Z96} 
affirme de plus que $Y$ est de type g\'en\'eral, ce que nous obtiendrons plus loin. 

\begin{demo}   
On se réduit aisément au cas où $S$ est  connexe, ce qu'on supposera dans la suite de la démonstration. 
Nous supposerons $S$ défini sur un corps de nombres $L$. 
Il n'est pas restrictif de supposer que l'image de $\gf{X}$ par $\rho$ 
est sans torsion. Si $g_\rho:X\To Y$ désigne la $\Gamma$-réduction associée à $\rho$,
 la représentation $\rho$ factorise  par $\gf{Y}$ d'après le lemme \ref{lemme sans torsion}. Considérons alors $H$ le sous-groupe de $\gf{Y}$ normalement engendré par les images des morphismes $\gf{Z}\To\gf{Y}$ induits par les applications holomorphes (avec $Z$ lisse et connexe) vérifiant $\rho(\gf{Z})=1$. Par construction, la représentation $\rho$ factorise par $\gf{Y}/H$, c'est à dire :
$$\rho\in Q:=\mathrm{Im}\left(M_B(\gf{Y}/H,S)\To M_B(Y,S)\right).$$
Notons au passage que $Q$ est définie sur $\bar{\QQ}$ et que les déformations de Simpson $([\rho_t])_{t\in\CC^*}$ sont également des éléments de $Q$.

Nous souhaitons factoriser toute la non rigidité de la représentation $\rho$. Pour cela, considérons tous les points $\sigma_\wp\in Q$ sur un corps local de la forme $L_\wp$ ainsi que les différentes composantes irréductibles $T_1,\dots,T_r$ de $Q$ qui passe par $[\rho]$. Le lemme \ref{factorisation KZ} fournit des applications holomorphes $s_{\sigma_\wp}$, $s_{T_i}$ et nous noterons $f:Y\To B$ la factorisation de Stein simultanée de toutes ces applications (un nombre fini d'entre elles suffit bien sûr à décrire $f$). Si $i:F\hookrightarrow Y$ désigne la fibre générale de $f$, la classe de conjugaison de la restriction de $\sigma$ à $\gf{F}$ est indépendante de $[\sigma]\in Q$ (avec $\sigma$ réductive). En particulier, la restriction de $\rho$ à $\gf{F}$ est conjuguée à une représentation définie sur un corps de nombres $L$ ; les déformations de Simpson étant constantes, $\rho_F$ est sous-jacente à une $\CC$-\textsc{vsh}. Enfin, pour tout place non archimédienne $\wp$ de $L$, $\rho_\wp:\gf{F}\To S(L_\wp)$ a 
une image précompacte, ce qui signifie exactement que la monodromie de $\rho_F$ est discrète.

Nous pouvons ainsi appliquer le lemme \ref{cas discret} : 
$\rho_F$ a un morphisme de Shafarevich $F\To Sh_{\rho_F}(F)$ dont la base est 
de type général. Par définition, la fibre générale de $F\To Sh_{\rho_F}(F)$ a
 une image triviale par $\rho$ ; mais $Y$ étant elle-même déjà la base de la 
$\Gamma$-réduction attachée \`a $\rho$, cela signifie exactement que la 
fibre générale de $F\To Sh_{\rho_F}(F)$ est un point et donc que $F$ est de type général. 
L'application $f$ est donc une fibration obtenue comme factorisation de Stein 
d'une application vers une variété de Kummer et sa fibre générale est de type général. 
Le théorème \ref{échange fibration} montre alors qu'il existe une application génériquement finie $r:Y'\To Y$ et une fibration en tores $g:Y'\To Z$ sur une variété de type général $Z$.

L'image de $r_*:\gf{Y'}\To\gf{Y}$ étant d'indice fini, 
la représentation $r^*\rho:\gf{Y'}\To S$ est encore Zariski dense, puisque $S$ est connexe. 
D'autre part, pour $z\in Z$ général, la fibre $Y'_z=g^{-1}(z)$ est un tore et 
la clôture de Zariski de $r^*\rho(\gf{Y'_z})$ est donc un sous-groupe abélien (puisque $\gf{Y'_z}$ l'est) 
et \emph{normal} de $S$ (puisque $\mathrm{Im}\left(\gf{Y'_z}\To\gf{Y'}\right)$ l'est). 
Le groupe $S$ étant supposé semisimple, $r^*\rho(\gf{Y'_z})$ est donc un groupe fini.
 La variété $Z$ étant la base d'une $\Gamma$-réduction associée à  $\rho$, nous en déduisons 
à nouveau que la fibre générale de $g$ est un point et que $Y'$ est de type général. \end{demo}

\subsection{Cas général}

Les résultats précédents s'étendent naturellement aux représentations réductives des groupes kählériens. Nous allons constater qu'elles sont construites à partir des variétés projectives et des tores complexes.\\

Soit donc $X$ une variété kählérienne compacte et $\rho:\gf{X}\To \gl_N(\bar{k})$ 
une représentation réductive définie sur un corps algébriquement clos de caractéristique zéro. 
L'adhérence de Zariski de $\rho(\gf{X})$ est donc un groupe algébrique réductif 
dont la composante neutre est  un sous-groupe d'indice fini. 
Le lemme suivant résulte donc de la description de la structure de ces groupes.

\begin{lem}\label{lemme réductif}
Il existe un revêtement étale fini $p_1:X'\To X$ tel que l'adhérence de Zariski de $G:=\rho(\gf{X'})$ soit un groupe réductif connexe. Par conséquent, il existe une isogénie canonique
$$i:G\To G^\sharp:=G/[G,G]\times G/Z(G)\simeq (\CC^*)^q\times S$$
où $S$ est semi-simple.
\end{lem}

\noindent Combiné avec les résultats de la partie précédente, ce lemme fournit la description suivante.

\begin{thm}\label{factorisation réductif}
Dans la situation décrite ci-dessus, il existe une modification d'un revêtement étale fini
 $\mu:{X}^\sharp\To X$, une isogénie $i:G\To G^\sharp$ comme en \ref{lemme réductif}, 
une fibration $s:{X}^\sharp\To Y$ sur une variété algébrique et une représentation réductive
$$\rho^\sharp:\gf{\alb({X}^\sharp)\times Y}\To G^\sharp$$
tels que $i\circ \mu^*\rho=(\alpha_{{X}^\sharp} \times s)^*\rho^\sharp$ 
(où $\alpha_{{X}^\sharp}:{X}^\sharp\To \alb({X}^\sharp)$ désigne l'application d'Albanese de ${X}^\sharp$).

\noindent Nous pouvons de plus imposer la condition $$\rho^\sharp\left(\gf{\alb({X}^\sharp)}\times \set{1}\right)\cap\rho^\sharp\left(\set{1}\times\gf{Y}\right)=\set{1}.$$
\end{thm}
\begin{demo}
Comme la partie abélienne de la représentation factorise 
nécessairement par l'application d'Albanese, il ne nous reste qu'à traiter la partie semisimple. 
Pour cela, appliquons le lemme \ref{lemme réductif} et considérons $\rho_s:\gf{X'}\To S$ 
la représentation obtenue en composant $\rho$ avec $i$ puis en projetant sur $S$. 
Le théorème \ref{factorisation presque simple} fournit une fibration méromorphe 
$X''\dashrightarrow Y_s$ sur un revêtement étale de $X'$ avec $Y_s$ algébrique. 
Si ${X}^\sharp\To Y$ désigne un modèle lisse de la factorisation de Stein 
de $X''\dashrightarrow Y_s$ ($Y$ est donc bien une variété algébrique lisse), 
la représentation $\gf{X^\sharp}\To S$ factorise par $Y$ (si on suppose que l'image de $\gf{X^\sharp}\To S$ 
est sans torsion, ce que nous ferons bien entendu) et cette fibration est bien celle recherchée.
\end{demo}

\subsection{Ubiquité des variations de structure de Hodge}

Comme rappelé dans la partie \ref{section Simpson}, une des conséquences
  de la théorie de C. Simpson est le phénomène d'ubiquité des variations de structures de Hodge. 
V. Koziarz a ainsi fait remarquer aux auteurs que la démonstration de Simpson ne se transpose pas
 au cas kählérien. En effet, les arguments de \cite{S92} reposent sur l'existence de $M_{Dol}(X,\gl_N)$ comme variété quasi-projective (si $X$ est algébrique), cet espace de module étant construit par les techniques de la Théorie Géométrique des Invariants. Les arguments développés ci-dessus permettent cependant de ramener le cas kählérien au cas projectif.

\begin{prop}\label{ubiquité kähler}
Soit $[\rho]\in M_B(X,\gl_N)(\CC)$ la classe de conjugaison d'une représentation réductive et soit $(\mathcal{E},\theta)$ le fibré de Higgs polystable associé. Pour $t\in \CC^*$, le fibré de Higgs polystable $(\mathcal{E},t\theta)$ est associé à une (classe de conjugaison de) représentation réductive $[\rho_t]\in M_B(X,\gl_N)(\CC)$. La limite
$$\lim_{t\to 0}[\rho_t]=[\rho_0]$$
existe dans $M_B(X,\gl_N)(\CC)$ et est sous-jacente à une $\CC$-\textsc{vsh}.
\end{prop}

\noindent Nous aurons besoin du lemme suivant.
\begin{lem}\label{propreté variété caractère}
Soit $\Gamma$ un groupe de type fini et $\Gamma_0\leq\Gamma$ un sous-groupe d'indice fini de $\Gamma$. L'application de restriction
$$\mathrm{Res}:M_B(\Gamma,\gl_N(\CC))\To M_B(\Gamma_0,\gl_N(\CC))$$
est propre (donc fini).
\end{lem}
\begin{demo}
Si $\Gamma=\gf{X}$ où $X$ est une variété kählérienne compacte,
 il s'agit d'une conséquence immédiate de \cite[Lemma 2.8]{S92}. 
Nous omettons l'argument dans le cas général puisque nous ne l'utiliserons pas. 
\end{demo}

\begin{demo}[du théorème \ref{ubiquité kähler}]
Considérons tout d'abord l'application $\mu:{X}^\sharp\To X$ 
et l'isogénie $i$ fournies par le théorème \ref{factorisation réductif} ; 
$\mu$ se décompose en une modification propre $b: {X}^\sharp\To X'$ et un revêtement
 étale fini $p:X'\To X$. De plus, ${X}^\sharp$ est l'espace total 
d'une fibration sur une variété algébrique lisse $Y$ et la représentation 
$\rho$ vérifie $i\circ \mu^*\rho=(\alpha_{{X}^\sharp} \times s)^*\rho^\sharp$ 
où $\rho^\sharp$ est une représentation de $\gf{\alb({X}^\sharp)\times Y}$ 
dans $G^\sharp$. En particulier $\rho^\sharp=\rho^\sharp_{ab}.\rho^\sharp_Y$
o\`u $\rho^\sharp_{ab}=p_{\alb({X}^\sharp)}^*\rho^\sharp|_{\gf{\alb({X}^\sharp}}$
et $\rho^\sharp_{Y}=p_Y^*\rho^\sharp|_{\gf{Y}}$, ces deux repr\'esentations \'etant vues comme
\`a valeurs dans deux sous groupes algébriques réductifs not\'es respectivement $G_1$ et $G_2$ qui 
de plus commutent entre eux. 
On a $[\rho^\sharp_t]=[\rho^\sharp_{ab,t}.\rho^\sharp_{Y,t}]$ o\`u $\rho^\sharp_{ab,t}$ 
(resp. $\rho^\sharp_{Y,t}$) est un repr\'esentant
de $[\rho^\sharp_{ab,t}]$ (resp. de $[\rho^\sharp_{Y,t}]$) \`a valeurs dans $G_1$ (resp. dans $G_2$). 
Comme $Y$ est projective lisse, nous concluons de la continuit\'e de $p^*_Y$
et du th\'eor\`eme d'ubiquit\'e pour $Y$ que la limite quand $t\to 0$
de $[\rho^\sharp_{Y,t}]$
existe. L'\'enonc\'e est \'el\'ementaire pour les repr\'esentations ab\'eliennes
et donc pour $ [\rho^\sharp_{ab,t}]$. Par suite
la limite quand $t\to 0$
de $[\rho^\sharp_t]$  existe. 
 Par image réciproque, la limite des déformations de Simpson de
 la représentation $i\circ \mu^*\rho$ (et donc $\mu^*\rho$ elle aussi\footnote{comme $i$ est une isogénie, le morphisme $M_B(\tilde{X},G)\To M_B(\tilde{X},G^\sharp)$ est fini.}) existe quand $t\to 0$. Comme $\mu^*\rho=b^*p^*\rho$ et 
comme $b_*$ est un isomorphisme entre les groupes fondamentaux, la conclusion reste valable pour $p^*\rho$. Nous pouvons finalement appliquer le lemme \ref{propreté variété caractère} pour constater que la limite $\lim_{t\to 0}[\rho_t]$ existe. Il est alors immédiat de vérifier que la classe d'isomorphisme du fibré de Higgs correspondant est fixée par l'action de $\CC^*$ et les arguments de \cite[section 4]{S92} montrent qu'un tel fibré correspond à une $\CC$-\textsc{vsh}.
\end{demo}

\subsection{Morphisme de Shafarevich}

Le théorème \ref{factorisation réductif} et l'existence du morphisme de Shafarevich démontrée dans \cite{E04}
dans le cas projectif impliquent:

\begin{prop}\label{repss} Soit $X$ une vari\'et\'e k\"ahl\'erienne compacte et $\rho: \gf{X} \to \Gamma=\rho(\gf{X})<\gl_N(\CC)$ une repr\'esentation 
lin\'eaire. 
Si $\Gamma<\gl_N(\CC)$ est d'adhérence de Zariski réductive, alors le morphisme de Shafarevich associé à $\rho$ existe.
\end{prop}

\begin{demo} L'existence du morphisme de Shafarevich est 
invariante par revêtement étale fini et modification biméromorphe. Par suite,  on peut supposer que l'adh\'erence de Zariski de $\Gamma$ est connexe et, en appliquant le
 théorème \ref{factorisation réductif}, on se r\'eduit au cas o\`u $X=X^{\sharp}$. On dispose alors
 d'une vari\'et\'e alg\'ebrique lisse $Y$ et d'une fibration holomorphe $s:X \to Y$ telle que la repr\'esentation $\rho$ factorise comme $( \alpha_X \times s)^*\rho^{\sharp}$
 o\`u $\rho^{\sharp}: \pi_1(\mathrm{Alb}(X) \times Y) \to \Gamma$ est une repr\'esentation. On note $\rho_{Y}$ la restriction de $\rho^{\sharp}$ \`a $\pi_1(Y)$. 
 Le morphisme de Shafarevich  $sh_{\rho_Y}: Y\to Sh_{\rho_Y}(Y)$ 
associé \`a $\rho_Y$ existe par \cite{E04}. On note $\rho^{ab}$ la restriction de $\rho^{\sharp}$ \`a $\pi_1(\mathrm{Alb}(X))$ et $sh_{\alpha_X^* \rho^{ab}}: X\to Sh_{\alpha_X^*\rho^{ab}}(X)$
le morphisme de Shafarevich correspondant dont la construction est donn\'ee \`a l'exemple \ref{exab}. 
On v\'erifie alors ais\'ement que le morphisme de Shafarevich  $X\to Sh_{\rho}(X)$ est la factorisation de Stein de  $sh_{\rho^{ab}} \times (sh_{\rho_Y}\circ  s)$.

\end{demo}

\section{Quotients résolubles linéaires des groupes k\"ahlériens}\label{th}

Nous étudions maintenant le cas d'un quotient linéaire \emph{résoluble} $\gf{X}\To\Gamma$ et établissons dans ce cas particulier
les résultats énoncés dans l'introduction.
 Les outils de base sont la théorie de Hodge, l'application d'Albanese et la structure de Hodge 
mixte sur la complétion nilpotente de $\gf{X}$, due à R. Hain \cite{Ha1,Ha2}.

\subsection{Morphisme d'Albanese}

\begin{prop}\label{qresol} Supposons $\Gamma$ résoluble (non nécessairement linéaire). 
Alors:
\begin{enumerate}[(1)]
\item Si $\Gamma$ n'est pas virtuellement nilpotent, il existe (à revêtement étale fini près) une 
application holomorphe surjective $f:X\to C$ sur une 
courbe projective de genre $g\geq 2$ qui factorise $\rho$.
\item Si $\Gamma$ n'est pas virtuellement abélien, il existe une application holomorphe 
$f:X\to Y$, o\`u $Y$ est une sous-variété de type général 
de dimension strictement positive d'une variété abélienne (quotient de $Alb(X)$).
\end{enumerate}
\end{prop}

\begin{demo}
L'assertion 1 est établie 
(à l'aide de l'invariant de Bieri-Neumann-Strebel) par T. Delzant dans \cite{Del06}.
 Le cas où $\Gamma$ est linéaire avait été traité (sans cet invariant) dans \cite{Ca01},
 basé sur les travaux classiques de Green-Lazarsfeld, Beauville, Simpson, Arapura-Nori.
L'assertion 2 (qui traite donc le cas où $\Gamma$ est virtuellement nilpotent) résulte de \cite{Ca95}.
\end{demo}

\begin{cor}\label{fsa} S'il n'existe pas de revêtement étale fini de $X$ admettant une application méromorphe surjective sur une variété de type général non-triviale, les quotients résolubles de $\pi_1(X)$ sont tous virtuellement abéliens.
\end{cor}

\begin{rem}\label{fspecspec}
\begin{enumerate}[1.]
\item L'hypothèse du corollaire équivaut \`a: l'application d'Albanese de tout revêtement étale fini de $X
$ est surjective et connexe.
\item Une variété spéciale au sens de \cite{Ca04}) satisfait cette hypothèse.
\end{enumerate}
\end{rem}

\subsection{Morphisme de Shafarevich}

\begin{thm} \label{reso}
 Soit $\rho:\pi_1(X)\to R$ un quotient résoluble linéaire. Alors, quitte à remplacer $X$ par un revêtement étale fini adéquat, $g_\rho:X\to G_\rho(X)$ coincide avec $g_{\rho^{ab}}:X\to G_{\rho^{ab}}(X)$, si $\rho^{ab}:\gf{X}\to R^{ab}$ est la représentation composée de $\rho$ avec le quotient $R\to R^{ab}$ de $R$ sur son abélianisé $R^{ab}$ (La fibration $g_{\rho^{ab}}$ est décrite \`a l'aide du morphisme d'Albanese dans l'exemple \ref{exab}).
\end{thm} 

\begin{rem} La démonstration qui suit ne permet pas d'aborder le cas résoluble non linéaire.
\end{rem}

\begin{demo} Quitte à remplacer $X$ par $X'$, revêtement étale fini, 
nous pouvons supposer que $D\Gamma:=[\Gamma,\Gamma]$, le groupe dérivé de $\Gamma$, 
est nilpotent\footnote{c'est ici que nous utilisons l'hypothèse de linéarité de $\Gamma$.}, noté $N$. Soit $R^{ab}:=\Gamma/D\Gamma$ l'abélianisé, 
et $\rho^{ab}:\gf{X}\to R^{ab}$ la représentation composée. 
Soit $g_{\rho/\rho^{ab}}:G_\rho(X)\to G_{\rho^{ab}}(X)$ la fibration introduite
 dans le lemme \ref{gred-comp}. Il s'agit de montrer que $g_{\rho/\rho^{ab}}$ est un isomorphisme. 
D'après ce lemme \ref{gred-comp}, la restriction de $g_{\rho^{ab}}$ \`a 
la fibre générique $X_z$ de $g_{\rho^{ab}}$ n'est autre que 
la $N$-réduction de $X_z$, si l'on note $\nu_z:\gf{X_z}_X\to N$ le morphisme 
(que l'on peut supposer surjectif, par le lemme \ref{gred-imrec}) obtenu par restriction 
de $\rho$ \`a $\gf{X_z}_X:=Im(\gf{X_z}\to \gf{X})$. Or la fonctorialité de la 
Structure de Hodge mixte sur la complétion de Mal\v{c}ev de $\gf{X_z}$
 combinée avec le caractère strict \cite{D} des 
morphismes de structures de Hodge mixte montre que l'application $\nu_z$ 
est d'image finie, d'après le lemme \ref{shm} ci-dessous. Ceci établit l'assertion.
\end{demo}

\begin{lem} \label{shm} Soit $Y$ une variété K\"ahlérienne compacte connexe et $Z$
un espace complexe compact k\"ahlérien connexe, possiblement non normal et non irréductible.
 Soit $f:Z\to Y$ une application holomorphe. Si l'application induite $f_*:H_1(Z,\Bbb Q)\to H_1(Y,\Bbb Q)$ est nulle, 
 alors, pour tout $k\geq 0$, $f_*(\gf{Z}/C^{k+1}\gf{Z})$ est fini, en notant $C^{k+1}G$ le $k$-ième terme de la suite centrale descendante d'un groupe $G$, définie par: $C^{k+1}G:=[G,C^kG], C^0G:=G$ . 
\end{lem} 
\begin{demo}
Il s'agit de \cite[cor. 5.2]{Ca98} si $Z$ est lisse. Le cas où $Z$ est normale en découle. Le cas général, 
d\^u  dans le cas projectif à Katzarkov \cite{Katnilp} et à S.~Leroy dans le cas k\"ahlérien (voir \cite{Cla08})
 résulte de \cite[Proposition 3.6, Remark 3.8]{EKPR}
appliquée
avec $M$ la représentation triviale et $G=GL(1)$.  
\end{demo}

Ce lemme permet en fait de montrer l'existence du morphisme de Shafarevich dans ce cas: 

\begin{cor} \label{resmshaf}Soit $\rho:\pi_1(X)\to R$ 
un quotient résoluble linéaire. Alors  le morphisme de Shafarevich $sh_\rho: X\to Sh_\rho(X)$ existe.  
Sur un revêtement étale fini 
adéquat de $X$, il coincide avec la factorisation de Stein de 
l'application composée $q\circ \alpha_X:X\to A_\rho$, o\`u $\alpha_X:X\to Alb(X)$ 
est le morphisme d'Albanese, et $q:Alb(X)\to A_X$ est 
le quotient par le plus grand des sous-tores compexes $T$ de $Alb(X)$ 
tels que $\gf{T}\subset K:=Ker(H_1(X,\Bbb Z)\to R^{ab})$. 
\end{cor}

\begin{demo} Sans perte de généralité,  nous pouvons supposer $\Gamma$ et $\Gamma^{ab}$ sans torsion. 
Considérons en effet la factorisation de 
Stein $g_\rho:X\to Y_\rho$ du morphisme $q\circ \alpha_X:X\to A_\rho$. 
Soit $f:Z\to X$ holomorphe, $Z$ compact connexe possiblement singulier et même non irréductible. 
Si $\rho\circ f_*(\gf{Z})$ est 
fini, alors il en est de même de $\rho^{ab}\circ f_*(\gf{Z})$, et $g_\rho\circ f(Z)$ est 
bien un point, puisque $g_\rho$ est le morphisme de Shafarevich de $\rho^{ab}$.
 Si $g_\rho\circ f(Z)$ est un point, alors $\rho^{ab}(f_*(\gf{Z}))$ est fini donc trivial. 
En particulier, $\rho (f_*(\gf{Z}))<N$, si $N:=D\Gamma<\Gamma$ est le groupe dérivé de $\Gamma$,
 qui est nilpotent. Du  lemme \ref{shm} précédent, appliqué \`a la composée: $q_\rho\circ \alpha_X\circ f:Z\to A_\rho$, nous déduisons que l'image de $\gf{Z}$ dans $N$, et donc dans $\Gamma$, est finie, ce qui établit que $g_\rho$ est le morphisme de Shafarevich de $\rho$.
\end{demo}

\subsection{Convexité holomorphe: cas linéaire résoluble géométrique}

Le cas linéaire résoluble se ramenant essentiellement au cas abélien, nous en déduisons l'énoncé de convexité holomorphe suivant.
\begin{cor} Soit $\rho:\pi_1(X)\to R$ un quotient résoluble linéaire. 
Supposons que $\rho^{ab}(\pi_1(X))$ soit  un quotient abélien géométrique. Alors, 
le revêtement $\widetilde{X_\rho}$ est holomorphiquement convexe.
\end{cor} 

\begin{demo} 
Il  résulte  du corollaire \ref{resmshaf} que l'application holomorphe naturelle
 $R(\widetilde{X_\rho})\to R(\widetilde{X_{\rho^{ab}}})$ 
est un rev\^etement topologique. Comme tout rev\^etement d'un espace de Stein est de Stein,
 l'exemple \ref{convhol} permet de conclure.
\end{demo}

\section{Convexité holomorphe des revêtements
linéaires des variétés kählériennes compactes}

Dans cette partie, nous donnons des conditions suffisantes sur $\rho$ garantissant la convexité holomorphe
de  $\widetilde{X_\rho}$.

\subsection{Revêtements réductifs des variétés kähériennes compactes}

Les résultats de \cite{E04} fournissent des conditions suffisantes de convexité holomorphe pour les revêtements réductifs des variétés projectives lisses et ce en terme de sous-ensembles constructibles absolus de $M_B(X,G)$ (avec $G$ un groupe réductif sur $\bar{\QQ}$). La généralisation de cette notion au cadre kählérien n'est \emph{a priori} pas évidente : celle-ci doit être fonctorielle pour les applications holomorphes entre variétés kählériennes compactes et pour les morphismes entre groupes réductifs (voir remarque \ref{remarque variété caractère}). Dans le cas abélien ($G=\gl_1$), les ensembles constructibles absolus doivent correspondre à des translatés de sous-tores par des points de torsion \cite{S93,Ca01}.

Avant de formuler une définition prenant en compte ces différentes exigences, fixons quelques notations. Considérons pour cela $T$ un tore complexe, $Y$ une variété projective lisse et $G/\bar{\QQ}$ un groupe algébrique réductif dont nous noterons $Z=Z(G)^\circ$ la composante neutre du centre. Si $\rho_T$ (resp. $\rho_Y$) désigne une représentation du groupe fondamental de $T$ (resp. de $Y$) à valeurs dans $Z$ (resp. dans $G$), le produit :
$$\rho_T\cdot\rho_Y(\gamma_T,\gamma_Y)=\rho_T(\gamma_T)\rho_Y(\gamma_Y)$$
définit bien une représentation de $\gf{T\times Y}$ à valeurs dans $G$. Après passage au quotient, nous obtenons un morphisme entre variétés des caractères :
\begin{align*}
M_B(T,Z)\times M_B(Y,G) &\To  M_B(T\times Y,G)\\
([\rho_T],[\rho_Y])&\longmapsto [\rho_T\cdot\rho_Y]=[\rho_Y]\cdot[\rho_Y]
\end{align*}
Pour finir, remarquons les points suivants : $M_B(T,Z)$ s'identifie au groupe algébrique $Z^{2\dim(T)}$ et le morphisme ci-dessus définit une action algébrique ; d'autre part, $M_B(Y,G)$ s'identifie au fermé $p_Y^*M_B(T\times Y,G)$ des représentations triviales en restriction à $\gf{T}$.

\begin{defi}\label{defi constructible absolu}
Soit $X$ une variété kählérienne compacte et $G$ un groupe algébrique linéaire réductif défini sur $\bar{\QQ}$.

Un ensemble constructible $M\subset M_B(X,G)(\bar{\QQ})$ est dit \emph{constructible absolu} si pour toute composante irréductible $\bar{M}$ de son adhérence de Zariski, il existe un revêtement étale $p:X'\To X$, une application méromorphe $s:X'\To T\times Y$ (avec $T$ un tore complexe et $Y$ une variété algébrique), un groupe algébrique réductif connexe $H\subset G$ et une isogénie $i:H\To H^\sharp$ tels que :
\begin{enumerate}
\item $p^*\bar{M}=j_*\bar{M}'$ où $j:H\To G$ est l'inclusion et $\bar{M}'$ un fermé de Zariski de $M_B(X',H)$,
\item $i_*\bar{M}'=\bar{s^*N}$ où $N\subset M_B(T\times Y, H^\sharp)$,
\item $N=M_B(T,Z(H^\sharp)^\circ)\cdot N_Y$ où $N_Y\subset M_B(Y,H^\sharp)$ est constructible absolu dans le sens de \cite{S93}.
\end{enumerate}
\end{defi}
\noindent Notons que $j_*$ étant un morphisme fini, il est en particulier fermé.

\begin{lem}\label{lemme defi équialente}
Si $X$ est algébrique, un ensemble constructible absolu dans le sens de \cite{S93} est un précisément un ensemble constructible absolu dans le sens de la définition \ref{defi constructible absolu}.
\end{lem}
\begin{demo}
L'implication directe est évidente en prenant pour $T$ un point, $Y=X$ et $p=s=\mathrm{id}_X$.

Réciproquement, si $X$ est algébrique, le tore $T$ peut être choisi algébrique et $M_B(T,Z(H^\sharp)^\circ)\cdot N_Y$ est alors constructible absolu au sens de Simpson (comme on le constate aisément). En particulier, $i_*\bar{M}'$ est constructible absolu au sens de Simpson. Le morphisme $i_*$ étant fini sur son image, $\bar{M}'$ est une composante irréductible de $(i_*)^{-1}(i_*(\bar{M}'))$ et donc constructible absolu au sens de Simpson. A nouveau, $p^*$ étant finie sur image (lemme \ref{propreté variété caractère}), le même argument montre que $\bar{M}$ est constructible absolu au sens de \cite{S93}.
\end{demo}

\begin{lem}\label{M_B est constructible absolu}
Si $X$ est une variété kählérienne compacte, $M_B(X,G)$
 et ses composantes irréductibles sont constructibles absolus.
\end{lem}
\begin{demo}
Soit $M$ une composante irréductible de $M_B(X,G)$ et $[\rho]$ un point générique de $M$. La représentation $\rho$ est alors conjuguée à une représentation Zariski dense $\rho_H$ à valeurs dans $H(\bar{k})$ où $\bar{k}$ est un corps algébriquement clos de degré de transcendance fini sur $\bar{\QQ}$ et $H\subset G$ est un $\bar{\QQ}$-sous-groupe. Nous pouvons alors appliquer le théorème \ref{factorisation réductif} à $\rho_H$. Prenons pour $N_Y$ la composante de $M_B(Y,H^\sharp)$ qui contient la restriction de $\rho_{H}^\sharp$ à $\gf{Y}$ ; de plus, quitte à considérer un revêtement étale fini de $T$ (et donc un revêtement étale fini de $X'$), nous pouvons supposer que ${\rho_{H}^\sharp}_{\vert \gf{T}}$ est à valeurs dans $Z(H^\sharp)^\circ$. Les conditions de la définition \ref{defi constructible absolu} sont alors satisfaites, en remarquant que $M_B(T,Z(H^\sharp)^\circ)\subset i_*M_B(T,Z(H)^\circ)$.
\end{demo}

Le lemme suivant n'a pas été explicité dans \cite{E04}. 

\begin{lem}\label{lemme points génériques}
Soit $X$ une variété kählérienne compacte et $M\subset M_B(X,G)(\bar{\QQ})$ constructible. Si $[\rho_1],\dots,[\rho_k]$ désignent les points génériques des composantes irréductibles de $\bar{M}$, alors le sous-groupe normal $H_M\lhd\gf{X}$ défini par :
$$H_M:=\bigcap_{[\rho]\in M}\ker(\rho)$$
vérifie
$$H_M=\bigcap_{i=1}^k \ker(\rho_i).$$
\end{lem}

\begin{demo}  
Cet énoncé n'est pas spécifique aux groupes kählériens. Plus généralement, on obtient immédiatement son analogue
 pour la variété des caractères de n'importe quel groupe de type fini $\Gamma$ à valeurs dans $GL_N$ 
en utilisant
que les fonctions $\{ \rho\mapsto \mathrm{Tr}( \rho(\gamma))\}_{\gamma\in \Gamma}$ engendrent l'anneau 
des fonctions de $M_B(\Gamma, GL_N)$ et la description de ses points  
sur un corps algébriquement clos de caractéristique nulle. On se ramène immédiatement à $G=GL_N$. 
\end{demo}

En particulier $H_M$ apparait comme le noyau d'une représentation 
linéaire réductive $\bar \rho=\oplus \rho_i$ à valeurs dans un corps de
caractéristique nulle de degré de transcendance fini sur $\mathbb{Q}$, 
donc aussi d'une représentation  linéaire réductive sur $\mathbb{C}$.

Nous pouvons maintenant conclure quant à la convexité holomorphe du revêtement obtenu en considérant toute une famille (constructible absolu) de représentations réductives.
\begin{thm}\label{convexité réductif}
Soit $X$ une variété kählérienne compacte et $M\subset M_B(X,G)(\bar{\QQ})$ constructible absolu. Le revêtement 
$$\widetilde{X_M}:=H_M \backslash \widetilde{X^u}$$
est alors holomorphiquement convexe.
\end{thm}

\begin{rem}
 La proposition et le lemme \ref{lemme points génériques} fournissent immédiatement que $\widetilde{X_M}$
a une fibration propre sur un espace complexe normal sans sous espaces compacts de dimension positive. 
S'il existe une représentation complexe d'adhérence de Zariski semisimple
 telle que $\widetilde{X_{\rho}}$ n'est pas holomorphiquement convexe est à notre connaissance une question ouverte.
\end{rem}

\begin{demo}
Nous remontons pas à pas le cours de la définition \ref{defi constructible absolu}. Si $X$ est projective, il s'agit de l'énoncé de \cite[Théorème 3]{E04}. Si $X=T\times Y$ et $M=N=M_B(T,Z(H^\sharp)^\circ)\cdot N_Y$, le revêtement correspondant à $M$ est :
$$\widetilde{X_M}=\widetilde{T^u}\times \widetilde{Y_{N_Y}}$$
et est donc holomorphiquement convexe comme produit d'un espace affine par une variété holomorphiquement convexe. Si $X=X'$ et si $M=s^*N$ (notations de la définition \ref{defi constructible absolu}), $\widetilde{X'_M}$ est alors propre sur $\widetilde{T^u}\times \widetilde{Y_{N_Y}}$. Examinons l'effet d'une isogénie et supposons que $X=X'$ et $M=M'$. Les sous-groupes définies par $M'$ et $i_*M'$ vérifient
$$H_{M'}\lhd H_{i_*M'}\lhd\gf{X}$$
et le lemme \ref{lemme points génériques} montre que l'indice de $H_{M'}$ dans $ H_{i_*M'}$ est alors fini (majoré par $\mathrm{deg}(i)^k$) et donc que
$$\widetilde{X'_{M'}}\To \widetilde{X'_{i_*M'}}$$ est un revêtement étale fini. Enfin, remarquons que
$$\widetilde{X'_{M'}}\To \widetilde{X_{M}}$$
est propre. Tout ceci montre bien la convexité holomorphe de $\widetilde{X_{M}}$.
\end{demo}

\subsection{Revêtements linéaires}
\begin{thm}\label{convexité linéaire}
Soit $X$ une variété kählérienne compacte, $G$ un groupe algébrique linéaire 
réductif défini sur $\QQ$ et $M=M_B(X,G)$ la variété des caractères associée.
 Considérons alors $\widetilde{H^{\infty}_M}$ le sous-groupe de $\gf{X}$ défini 
comme l'intersection des noyaux de \emph{toutes} les représentations $\gf{X}\To G(A)$ 
où $A$ est une $\CC$-algèbre arbitraire. Le revêtement galoisien de groupe $\widetilde{H^{\infty}_M}$
$$\widetilde{X^{\infty}_M}:=\widetilde{X^u}/\widetilde{H^{\infty}_M}$$
est alors holomorphiquement convexe.
\end{thm}
\begin{rem}
 On a $\widetilde{H^{\infty}_M}=\ker(\rho_A^{taut})$ o\`u $A$ est l'anneau des fonctions r\'eguli\`eres sur
$R_B(\gf{X}, G)$ et $\rho^{taut}:\pi_1(X)\to G(A)$ est la repr\'esentation tautologique.
\end{rem}

\begin{cor}\label{corollaire linéaire}
Si $X$ est une variété kählérienne compacte dont le groupe fondamental est linéaire, son revêtement universel $\widetilde{X^u}$ est holomorphiquement convexe.
\end{cor}

\begin{demo}[du théorème \ref{convexité linéaire}]
La preuve donnée dans le cas projectif \cite{EKPR} s'applique modulo l'adaptation du
 lemme 5.1 de cette r\'ef\'erence au cas k\"ahlérien. Cette adaptation est ais\'ee, voir \cite{E13}. 
\end{demo}

\section{Structure des variétés de Shafarevich}

L'objectif de cette section est de démontrer 
un théorème de structure pour la variété de Shafarevich associée à une représentation linéaire $\rho$ en combinant les énoncés similaires 
dans les cas particuliers antithétiques et 
complémentaires où $\Gamma$ est (d'adhérence de Zariski) 
semi-simple d'une part, résoluble d'autre part (La réduction du cas général à ces deux cas particuliers résultant de la décomposition de Levi-Mal\v{c}ev).

\subsection{Variétés de Shafarevich semisimples}

Nous donnons ici notre d\'emonstration alternative de la derni\`ere partie du th\'eor\`eme principal de
\cite{Z96} en montrant que la base du morphisme de Shafarevich
 d'une représentation semi-simple est de type général. Pour cela, nous allons 
utiliser une construction élaborée dans \cite[\S 5.3]{E04}, celle du I-morphisme de Shafarevich. Rappelons en brièvement les grandes lignes : si $M$ est un sous-ensemble constructibe absolu de $M_B(X,G)(\bar{\QQ})$ (avec $G$ un groupe réductif et $X$ projective lisse), il est possible de lui associer
$$sh^I_M:X\to sh^I_M(X)$$
qui est obtenu comme factorisation de Stein d'un morphisme vers une certaine variété de Kummer comme dans la preuve du
 théorème \ref{factorisation presque simple}. Plus précisément, il s'agit de la réduction de Katzarkov-Zuo
 associée \`a un produit de représentations de $\pi_1(X)$ vers les points de $G$ sur divers corps locaux 
auquelles on peut attacher une application harmonique équivariante $h_M:\widetilde{X_M} \to \Delta$ o\`u
$\Delta$ un produit d'immeubles de Bruhat-Tits. Par \cite[Prop 5.4.6]{E04}, 
le feuilletage défini par le noyau de la différentielle complexifiée de $h_M$ 
coincide avec celui qu'induit la fibration $sh^I_M$. En particulier son corang\footnote{L'argument de \cite{Z96} repose sur l'énoncé plus fort que 
cette propri\'et\'e est satisfaite par la factorisation de Katzarkov-Zuo 
d'une  repr\'esentation 
Zariski dense dans un corps local non archim\'edien et en derni\`ere analyse
 sur \cite[Theorem 4.2.3]{Z99}.  Une difficulté dans l'argument de \cite[p. 72]{Z99} a  conduit \cite{E04}
a en donner la pr\'esente version affaiblie et a motivé le développement de notre approche alternative
à \cite{Z96}.} au point général est
 $r:=\dim (sh^I_M(X))$. 
Cet \'enonc\'e permet d'établir le lemme suivant:

\begin{lem}\label{type general}
Soit $X$ une variété kählérienne compacte et $M\subset M_B(X,G)(\bar{\QQ})$ constructible absolu. Soit 
$$\widetilde{X_M}:=H_M \backslash \widetilde{X^u}\to \widetilde{S_M(X)}$$
la réduction de Remmert de Stein, $\Gamma_M= \pi_1(X)/H_M$ et $X\to sh_M(X)=\Gamma_M \backslash \widetilde{S_M(X)}$ 
le morphisme de Shafarevich.

Si $\kappa(sh_M(X))=0$, ces représentations génériques sont d'image finie et $\Gamma_M$ est fini. 

Si $\Gamma_M$ est sans torsion et si les points génériques des composantes connexes de $M$ sont des caractères
de représentations 
d'adhérence de Zariski semi simple, $sh_M(X)$ est de type général.
\end{lem}

\begin{demo} Nous pouvons supposer que $M$ est irréductible, que $sh_M$ est biméromorphe (quitte à remplacer $X$ par la base du morphisme de Shafarevich correspondant) et que $G$ presque simple. Le théorème \ref{factorisation presque simple} nous réduit alors au cas où $X$ est projectif.

Pour montrer le résultat annoncé, il nous suffit d'établir les deux points suivants :
\begin{enumerate}[(I)]
\item la dimension de Kodaira de $X$ est positive\footnote{Notons que ce point résulterait de la conjecture d'abondance.} : $\kappa(X)\ge0$.
\item si $\kappa(X)=0$ , la représentation correspondant au point générique de $M$ est d'image finie.
\end{enumerate}
Voyons tout d'abord comment ceci permet de conclure que $X$ est de type général. Si $\kappa(X)\ge0$, il nous suffit de montrer que les fibres générales de la fibration d'Iitaka de $X$ sont des points. Si $F$ désigne une telle fibre, elle vérifie par définition $\kappa(F)=0$ et son image dans $G$ est un sous-groupe normal (car $F$ est la fibre générale d'une application holomorphe). Comme $G$ est presque simple, le deuxième point ci-dessus nous assure que l'image de $\gf{F}$ dans $G$ est finie. Comme nous avons supposé que $X=sh_M(X)$, cela implique que $F$ est réduit à un point et $X$ est bien de type général.\\

Pour établir les points ci-dessus, utilisons une construction de \cite[Paragraphe 1.3.4]{E04}. Si $\beta$ est une forme $r$-linéaire alternée
complexe sur l'appartement de $\Delta$, la quantité $<\prod_{w\in W} w^* \beta ,\Lambda^r\partial h_M>$ 
définit un tenseur holomorphe  $\tau \in H^0(X, S^N \Omega^r_X)$ qui est non nul si $\beta$ est choisie génériquement.

 Par construction $\tau$ est induit par une forme pluricanonique holomorphe $\tau^I$ définie sur un ouvert 
$U\subset sh^I_M(X)$. En choisissant une multisection de $sh^I_M$ nous voyons qu'il existe une 
application holomorphe $\psi:X'\to sh^I_M(X)$, $X'$ étant projective lisse, telle que 
$\psi^* \tau^I$ se prolonge \`a une forme pluricanonique régulière sur $X'$. Ceci signifie précisément 
que $\tau^I$ est une forme pluricanonique sur la base orbifolde de $sh^I_M$ \cite{Ca04}, 
c'est-à-dire $\kappa\left(sh^I_M(X),\Delta(sh^I_M)\right)\ge0$. Remarquons également que 
les fibres $F$ de $sh^I_M$ sont de type général (par construction elles ont une application 
de périodes génériquement immersive, comme dans la démonstration du théorème \ref{factorisation presque simple}). 
Le  point (I):
$$\kappa(X)\ge \kappa(F)+\kappa\left(sh^I_M(X),\Delta(sh^I_M)\right)\ge0.$$
résulte alors de la variante orbifolde de \cite{Ko}:

\begin{lem}\label{addorb}
 Soit $f:X\dasharrow Y$ une  fibration rationnelle de variétés projectives complexes 
dont la fibre générale fibre est de type général. Alors,
$\kappa(X)\ge \dim(X_y)+\kappa(Y,f)$ o\`u $\kappa(Y,f)=\kappa(Y,\Delta(\hat{f}))$ et
$\hat{f}$ d\'esigne un mod\`ele net pr\'eparé et admissible de $f$ (voir \cite[Section 1.3]{Ca04}).
\end{lem}

\begin{demo}
Il s'agit de reprendre les arguments de \cite{V83,Ko} en y incorporant les considérations développées dans \cite{Ca04}. En effet, les arguments de \cite[\S 7]{V83} montrent qu'il suffit de traiter le cas d'une fibration dont la variation est maximale, une fois que l'on y a remplacé le théorème III par le résultat de semi-positivité orbifolde \cite[Th. 4.11]{Ca04}. Si $f$ est à variation maximale, il nous faut alors établir l'égalité :
\begin{equation}\label{image directe big}
\kappa\left(Y,\mathrm{det}f_*(mK_{X/(Y,\Delta(f))})\right)=\dim(Y)
\end{equation}
pour $m>>1$ (toujours d'après les arguments de Viehweg). En fait, dans cette situation, le faisceau $f_*(mK_{X/(Y,\Delta(f))})$ est lui même \emph{big}\footnote{nous renvoyons à \cite{V83} pour les notions de faible positivité et de \emph{bigness}.} pour $m$ assez grand. En effet, d'après \cite[Prop. 4.15]{Ca04}, il existe un changement de base (avec $v$) fini
$$\xymatrix{X'\ar[d]_u \ar[r]^{g} & Y'\ar[d]^v\\
X\ar[r]^f & Y
}$$
et une injection de faisceau
$$g_*(mK_{X'/Y'})\hookrightarrow v^*f_*(mK_{X/(Y,\Delta(f)}+B)$$
(avec $B$ un diviseur sur $X$ qui est $f$-exceptionnel et qui ne perturbe donc pas les espaces de sections globales puisque $f$ est supposée nette). Or, le résultat principal de \cite{Ko} établit le caractère \emph{big} du faisceau $g_*(mK_{X'/Y'})$ pour $m$ assez grand ; les faisceaux considérés étant de même rang, nous en déduisons que $v^*f_*(mK_{X/(Y,\Delta(f)}+B)$ ainsi que $f_*(mK_{X/(Y,\Delta(f)})$ sont également \emph{big} et que l'égalité (\ref{image directe big}) est bien vérifiée.
\end{demo}

Pour le point (II), nous reprendrons essentiellement les arguments de \cite{Z96}.
On suppose maintenant $\kappa(X)=0$. Le raisonnement mené ci-dessus
 montre que dans ce cas $sh^I_M$ est une application birationnelle et le noyau de la différentielle complexifiée de $h_M:\tilde{X}\To \Delta_M$ est trivial au point général de $\tilde{X}$. Dans ce cas, la construction ci-dessus fournit une forme pluricanonique qui s'annule sur le diviseur de ramification du revêtement spectral $X^s_M$ associé à $M$ (voir \cite{Z96,E04}) et nous en déduisons l'égalité : $\kappa(X^s_M)=\kappa(X)=0$. D'autre part, le revêtement spectral vient également avec une application vers une variété abélienne $\alpha_M:X^s_M\To \alb_M$ ; comme celle-ci provient essentiellement de l'application pluriharmonique $h_M$, l'hypothèse sur le noyau de la différentielle de $h_M$ montre que $\alpha_M$ est génériquement finie. Les résultats de \cite{K81} montrent que $X^s_M$ est alors birationnelle à une variété abélienne. Cette dernière remarque entraîne que l'image de la représentation correspondant au point générique de $M$ est virtuellement abélienne : elle est donc finie d'après \cite{S93}.
\end{demo}

\noindent Nous retrouvons ainsi le dernier point du théorème principal 
de \cite{Z96}.  

\begin{thm}\label{repss} 
Si $\Gamma<\gl_N(\CC)$ est d'adhérence de Zariski semi-simple et $\Gamma$ sans torsion, $Sh_\rho(X)$ est une variété de type général. 
\end{thm}

\begin{demo}
Si  $M$ le plus petit fermé constructible absolu de $M_B(X,\bar\Gamma)$ contenant $\rho$,
$M$ est irréducible. Consid\'erons  $\rho_M$ un représentant de son
 point générique. Puisque $\ker(\rho_M)<\ker(\rho)$,
 on a un morphisme $\tau: Sh_M(X)=Sh_{\rho_M}(X)\to Sh_{\rho}(X)$.  Soit $i:Z\to X$  la fibre générale de 
$X\to Sh_{\rho(X)}$ de sorte  que $i^*\rho$ est la représentation triviale (en particulier, $Z$ est  lisse).  La  classe de conjugaison
$[\mathds{1}]$
 de la représentation  triviale
est fermée absolue dans $M_B(Z, \bar \Gamma)$, $M\subset (i^*)^{-1}[\mathds{1}]$ ce qui signifie que $\rho_M(\pi_1(Z))=\{1\}$. 
Donc $Z$ est contractée dans $Sh_M(X)$ ce qui implique que $\tau$ est biméromorphe. Le lemme  \ref{type general}
permet alors de conclure. 
\end{demo}
\begin{rem}
Dans ce raisonnement, on peut prendre pour $Z$ un revêtement étale de la
 résolution des singularités d'une composante d'une fibre de  $\tau$.
 Il en résulte que $\tau$ est un isomorphisme. Néanmoins, nous ignorons s'il est possible 
que $\widetilde{X_{\rho}}$ ne soit pas holomorphiquement convexe, c'est \`a dire
que $R(\widetilde{X_{\rho}})$ soit de Stein. Ce que nous obtenons est que 
$R(\widetilde{X_{\rho_M}}) \to R(\widetilde{X_{\rho}})$ est un rev\^etement topologique
par un espace de Stein. 
\end{rem}

\subsection{Cas linéaire général}

Nous considérons maintenant le cas général $\rho:\gf{X}\to \Gamma<\gl_N(\Bbb C)$. Nous noterons $G:=\overline{\Gamma}$ l'adhérence de Zariski de $\Gamma$. Quitte à remplacer $X$ par un revêtement étale fini adéquat, nous supposerons que les groupes algébriques intervenant dans la suite sont connexes et que les représentations linéaires déduites de $\rho$, de ses quotients ou sous-représentations se factorisent par les $\Gamma$-réductions associées.

Soit donc $\bar R$ le radical résoluble de $G$ et $\lambda: G\to S:=G/R$ son quotient semi-simple. Par Levi-Mal\v{c}ev, $G$ est produit semi-direct de $S$ par $\bar R$. Nous désignerons par $R:=\Gamma\cap \bar R$, et $s:\Gamma\to \Gamma/R$ le quotient naturel déduit de $\lambda$. Donc $R$ est résoluble, et $\sigma:=s\circ \rho: \gf{X}\to \Gamma/R$ est semi-simple.

Nous avons donc (par le lemme \ref{gred-comp}) des réductions $g_\rho:X\to G_\rho(X)$ et $g_\sigma:X\to G_\sigma(X)$ associées à $\rho$ et $\sigma$ et une factorisation $g_{\rho/\sigma}:G_\rho(X)\to G_\sigma(X)$ telle que $g_\sigma=g_{\rho/\sigma}\circ g_{\rho}$. La restriction de $g_\rho$ à la fibre générique $X_\sigma$ de $g_\sigma$ n'est autre que la $R$-réduction de $g_{\rho_\sigma}$. En effet (voir lemme \ref{gred-comp}), la restriction, notée $\rho_\sigma$, de $\rho$ à $\gf{X_\sigma}$ a pour image $R$. 

En combinant le théorème \ref{repss} avec le corollaire \ref{resmshaf}, nous pouvons établir le :

\begin{thm}\label{t-lingen} 
Soit $\rho:\gf{X}\to \Gamma<\gl_N(\CC)$ comme ci-dessus. 
Alors (quitte à remplacer $X$ par une modification propre d'un revêtement étale fini adéquat), 
la variété $Sh_\rho(X)$ est biméromorphe à l'espace total d'une fibration lisse $\tau:Sh_\rho(X)\to S_\rho(X)$ 
en tores complexes sur une variété algébrique $S_\rho(X)$ de type général.
\end{thm}

\begin{demo} 
Notons $Sh_{\sigma}(X):=Z$. Posons $H=R/[R,R]$, $\rho':\pi_1(X) \to \Gamma/ [R,R]$ et appliquons le lemme \ref{gred-ab}. 
La $\rho'$-réduction $sh_{\rho'}:X\to Sh_\rho(X):=Y$ est 
biméromorphe \`a la factorisation de Stein d'un 
 $q_\rho\circ \alpha_{X/Z}:X\to A_\rho(X/Z):=A_\rho$, 
o\`u $\tau: A_\rho(X/Z)\to Z$ est 
une fibration presque holomorphe dont les fibres générales sont des tores complexes. 
Il résulte du lemme \ref{gred-comp} et du théorème \ref{reso} que cette $\rho'$ réduction
est la $\rho$-réduction de $X$. 

Nous avons, pour $z\in Z$ général une application 
de Ueno-Kawamata unique: $u_z: Y_z\to W_z$, qui est un fibré 
principal de groupe $K_z$, un revêtement étale fini d'un sous-tore de $A_\rho(X/Z)_z$,  
avec $W_z$ de type général. Par la compacité des composantes irréductibles 
de l'espace des cycles de Chow-Barlet de $Y$, nous en déduisons l'existence d'une 
application presque-holomorphe $u:Y\to W$ au-dessus de $Z$, 
holomorphe au-dessus de l'ouvert de lissité de $sh_\sigma: X\to Z$, 
et dont la restriction \`a $Y_z$ coïncide avec $u_z:Y_z\to W_z$, fibre de $W$ au-dessus de $z$.  Nous avons alors une application d'Albanese relative naturelle: $\alpha_{u}:\alb(Y/Z)\to \alb(W/Z)$ au-dessus de $Z$, et un morphisme 
injectif: $\alpha_{W/Z}:W\to \alb(W/Z)$ au-dessus de $Z$ tels 
que $\alpha_u\circ\alpha_{Y/Z}=\alpha_{W/Z}\circ\alpha_{Y/Z}:Y\to \alb(W/Z)$. 
L'application $u:Y\to W$ est surjective et une submersion \`a fibres 
des tores complexes au-dessus de l'ouvert $U\subset Z$ de lissité de $sh_\sigma$. 

La variété $S_\rho(X):=W$ est bien de type général, 
puisque fibrée sur $Z=Sh_\sigma$ de type général, 
et \`a fibres génériques $W_z$ de type général.

Pour conclure que la fibration obtenue $\tau:Sh_\rho(X)\to S_\rho(X)$ est biméromorphe
à une fibration lisse, il suffit de remarquer que le groupe fondamental de la fibre générale de $\tau$ s'injecte dans celui de $Sh_\rho(X)$ (ceci vient du fait que $Sh_\rho(X)$ est la base d'un morphisme de Shafarevich). Nous pouvons alors appliquer l'énoncé \cite[Th. 7.8]{Nak} qui nous fournit la conclusion souhaitée.
\end{demo}

\providecommand{\bysame}{\leavevmode\hbox to3em{\hrulefill}\thinspace}
\providecommand{\MR}{\relax\ifhmode\unskip\space\fi MR }
\providecommand{\MRhref}[2]{%
  \href{http://www.ams.org/mathscinet-getitem?mr=#1}{#2}
}
\providecommand{\href}[2]{#2}

\end{document}